\documentclass[preprint,12pt,authoryear]{elsarticle}

\usepackage[margin=0.9in]{geometry}
\usepackage[utf8]{inputenc}
\usepackage[english]{babel}
\usepackage{amsmath,amssymb,amsthm}
\allowdisplaybreaks % Place this in your preamble
\usepackage{graphicx}
\usepackage{float}
\usepackage{makecell}
\usepackage{algorithm}
\usepackage{algpseudocode}
\usepackage{booktabs}
\usepackage{multirow}
\usepackage{longtable}
\usepackage{tabularx}
\usepackage[hidelinks]{hyperref}

\biboptions{authoryear}

\newtheorem{proposition}{Proposition}
\theoremstyle{remark}
\newtheorem{remark}{Remark}

\journal{Preprint}

\begin{document}

\begin{frontmatter}

\title{Bilevel subsidy-enabled mobility hub network design with perturbed utility coalitional choice-based assignment}

\author[nyu]{Hai Yang}
\author[nyu]{Joseph J. Y. Chow\corref{cor1}}
\cortext[cor1]{Corresponding author}
\ead{joseph.chow@nyu.edu}

\address[nyu]{C2SMARTER University Transportation Center, Department of Civil \& Urban Engineering, New York University Tandon School of Engineering, 6 MetroTech Center, Brooklyn, NY 11201, USA}

\begin{abstract}
Urban mobility is undergoing rapid transformation with the emergence of new services. Mobility hubs (MHs) have been proposed as physical-digital convergence points, offering a range of public and private mobility options in close proximity. By supporting Mobility-as-a-Service, these hubs can serve as focal points where travel decisions intersect with operator strategies. We develop a bilevel MH platform design model that treats MHs as control levers. The upper level (platform) maximizes revenue or flow by setting subsidies to incentivize last-mile operators; the lower level captures joint traveler–operator decisions with a link-based Perturbed Utility Route Choice (PURC) assignment, yielding a strictly convex quadratic program. We reformulate the bilevel problem to a single-level program via the KKT conditions of the lower level and solve it with a gap–penalty method and an iterative warm-start scheme that exploits the computationally cheap lower-level problem. Numerical experiments on a toy network and a Long Island Rail Road (LIRR) case (244 nodes, 469 links, 78 ODs) show that the method attains sub–1\% optimality gaps in minutes. In the base LIRR case, the model allows policymakers to quantify the social surplus value of a MH, or the value of enabling subsidy or regulating the microtransit operator's pricing. Comparing link-based subsidies to hub-based subsidies, the latter is computationally more expensive but offers an easier mechanism for comparison and control.
\end{abstract}


\end{frontmatter}

\section{Introduction}

Urban mobility has significantly evolved due to advances in mobile technology and digital platforms. Yet, transportation networks continue to struggle with persistent challenges including congestion, emissions, and unequal accessibility. Despite the proliferation of new mobility services like ride-hailing, car-sharing, and micromobility options, there remain challenges to fully integrate these systems to fundamentally resolve urban mobility issues. Travelers still face inconvenience in seamlessly coordinating multiple modes at transfer locations. Large proportions of commuters may still prefer driving alone if the connections to major public transit stops are inadequate or unreliable.

Emergent cyberphysical mobility hubs (MHs) are a viable response to these challenges. MHs are defined as "multimodal transport nodes that facilitate intermodal transfers by providing different mobility options in close proximity" \citep{miramontes2017impacts}. While "transit hubs" have existed for decades as a means to integrate land use and public transport as a type of "transit-oriented development" (see \cite{weustenenk2023towards}), they differ from MHs in the literature in the last decade, which feature higher levels of both physical and digital integration of mobility services \citep{arias2023exploring}, particularly services that feature real-time information. These MHs often involve multiple operators (or even non-mobility-related goods and service providers, such as fueling), necessitating distinct incentives and cooperation mechanisms enabled by integrated digital platforms supporting Mobility-as-a-Service (MaaS). MaaS platforms streamline users' abilities to plan, book, and pay for various mobility services through a unified digital interface and distribute those revenues between operators, thereby aligning physical infrastructure and digital coordination. As such, MHs form critical components of cyberphysical MaaS ecosystems, acting as strategic nodes for interactions among public entities, private mobility operators, and users.

The deployment of MHs alongside MaaS schemes has been studied globally over the past decade. A notable example is the car-sharing program in Bremen, Germany \citep{KarbaumerWeltring2025}. Car-sharing stations known as “mobil.punkte” are established at easily accessible locations near major public transit stops, with supporting infrastructure such as bike racks and app-based booking systems. Sustainability-focused mobility hub designs have gained more attention with the maturation of electric vehicle (EV) technologies. The eHUBS project, supported by the European Regional Development Fund, investigated the deployment of intermodal hubs focusing on electric-powered modes in six pilot cities across North-West Europe \citep{eHUBS2025}. For a comprehensive review of mobility hub-focused pilot studies worldwide, readers can refer to \cite{arnold2023mobility}.

Designing systems incorporating MHs is no trivial task. \cite{grigolon2025willingness} reveals that travelers are willing to pay extra to use MHs when multiple modes are easily accessible through a well-designed physical location supported by a digital platform. Existing studies predominantly focus only on their role as physical transfer points (i.e. transit hubs), where the main objective is to optimally choose the hub locations that enhance accessibility \citep{petrovic2019location, frank2021improving, aydin2022planning}. However, these studies often overlook the potential benefits of MHs on operator strategies. For example, MHs provide a location for stationing on-demand vehicles, fueling or charging electric vehicles, and provide a means for two operators to have controlled cost transfers (i.e. only subsidize trips that start or end at MHs). For example, an integration between Long Island Railroad (LIRR) and Uber may be costly if LIRR subsidized all Uber trips to all stations, but setting boundaries for specific stations as designated MHs through bundled tickets purchased for Uber and LIRR can produce partnerships in a more cost-effective manner that can suit both parties. The problem becomes even more complex should other mobility provides wish to participate, i.e. how should LIRR consider a common platform linking with both Uber and Lyft, how should subsidies to these operators be designed, and which stations should be selected as MHs?

Effective service integration requires joint operational planning among operators to improve traveler experiences. Recent studies, such as \cite{xanthopoulos2024optimization}, partially address this gap by considering passenger preferences in determining optimal locations and capacities of mobility hubs. Nevertheless, these works still do not fully explore how such hub characteristics as location, capacity, pricing structures, and subsidies can directly shape operators’ key decisions, including fleet allocation, service area planning, and profitability. In reality, MHs can serve as powerful control levers within public cyberphysical platforms, critically influencing both demand-side traveler choices and supply-side operator strategies. This dual influence remains largely underexplored in the current literature.

Recent advancements in MaaS modeling provide a feasible way in addressing the complex interactions between travelers and operators. Bilevel frameworks and many-to-many stable matching games \citep{liu2024demand, yao2024design, liu2024modeling} have effectively captured market dynamics between travelers and operators. However, these models do not explicitly incorporate mobility hubs as integral network gateways. This gap highlights the need to systematically integrate mobility hub considerations into MaaS modeling frameworks, particularly in contexts involving complex interactions between multiple public and private operators, such as providing first-mile and last-mile services to commuter rail stations. In addition, high computational complexity associated with the MaaS models prevents larger scale deployment, which is often associated with MH design scenarios that must link both neighborhood subnetworks with larger regional subnetworks.

We fill this gap in the literature by proposing a mathematical model that optimizes the capacity allocation of shared multimodal mobility hubs while maximizing travel utility, accounting for multimodal trips and the cost transfers and subsidies between operators. In this context, operators participating in the designated MHs are part of a MaaS platform, which can be operated by a lead fixed route transit operator. To address computational complexity issues and scalability limitations inherent in existing multimodal flow assignment methods, we introduce the Perturbed Utility Route Choice (PURC)-based assignment game, based on the route choice framework from \cite{fosgerau2022perturbed}. Unlike path-enumeration approaches, the PURC framework uses link-level utilities to model dispersed multimodal route choices efficiently. This approach facilitates scalable and realistic MaaS market analyses involving mobility hubs. A case study for the neighborhoods along three LIRR stations in Suffolk County, NY, illustrates the application of the proposed mobility hub design model on a real-world scale. In addition, the case study provides valuable insights to how a MH can act as leverage in controlling the amount of subsidies used to improve transit station usage.

The paper is structured as follows. We present a literature review of the studies focusing on MH and multi-modal mobility network design games in the following section. We then explain the details of the methodology. Afterwards, we present the LIRR based case study to illustrate the model application and discuss the nuance of the service strategies involved in establishing mobility hubs. Finally, we conclude the study and discuss how the work can be extended.

\section{Literature Review}

Multimodal systems and mobility hub related studies cover several knowledge fields including network pricing modeling, facility location modeling, discrete choice modeling, spatial analysis, social equity evaluations, game-theory based stable matching, and many more. We do not intend to provide an exhaustive list of literature that covers all aspects of the related studies. Instead, we mainly focus on methodological works that explore the interactions among different participants when designing cooperative multimodal services (i.e. MaaS) or MH systems.

\subsection{MaaS platform models}
MaaS platforms integrate various transportation modes into bundles that enable more seamless trips involving multiple service operators. A MaaS platform allows users to purchase trips consisting of combinations of mobility options in a single transaction. The pricing and service schemes have attracted considerable attention in recent years \citep{van2022business}. \cite{bertsimas2020joint} proposed a framework that jointly optimizes schedule frequencies and service pricing in a multi-modal transit network. Given a set of multimodal trip paths, service providers optimize their service frequencies to minimize system-wide wait time under budget constraints, while passengers select their preferred routes based on a nested logit model. \cite{horcher2020maas} compared how different MaaS subscription policies affect system performance, specifically regarding net welfare, externalities, and operator objectives. The study involves three levels of decision making: car ownership, MaaS package selection, and daily mode choice. Each level is modeled using utility maximization principles with nonlinear externalities. They found that certain pricing policies can lead to unexpected negative outcomes by encouraging the overconsumption of services.

Alternatively, two-sided market matching mechanism is adopted to model MaaS service design problems. \cite{djavadian2017agent} pioneered the two-sided market mechanism into MaaS modeling. \cite{chen2023stability} studied the convergence properties of these two-sided markets with a focus on ride-hail matching. \cite{rasulkhani2019route} and \cite{pantelidis2020many} expanded the stable matching approach with transferable utilities to model the interaction of fixed-line service providers and travelers as a two-sided market, first as a many-to-one assignment game and then as a many-to-many assignment game. Travelers select their desired paths, which may involve multiple operators, while service operators determine whether to operate a service link and set service prices. \cite{liu2024demand} further extended this framework to incorporate mobility-on-demand (MOD) services with congested links (i.e. macroscopic matching impedance functions) for accessing MOD services. An exact solution algorithm was proposed to identify stable outcomes and determine the level of subsidy required from the platform to the operators to stabilize empty cores, highlighting the additional resources needed to operate MaaS platforms. Building on the many-to-many matching framework, \cite{yao2024design} proposed a different design in which the MaaS platform acts as an intermediary which purchases capacity from service providers and offers service bundles to travelers with OD-based pricing. Additionally, non-MaaS travelers are considered within the congested network.

Bilevel structures are also frequently used to formulate MaaS platform design problems. The upper level typically represents decisions made by the service platform or providers, while the lower level captures the decisions of other participants, such as travelers. \cite{xi2024single} proposed a single-leader multi-follower game (SLMFG) in which the MaaS platform is the leader, and travelers and service providers are followers. The platform determines service bundles and pricing, while travelers and service providers choose their participation to optimize their respective objectives. Building on this work, \cite{xi2024strategizing} expanded the framework into a multi-leader multi-follower game (MLMFG) that considers competition among platforms in an electric MaaS (E-MaaS) ecosystem. \cite{huang2024non} used an SLMFG to model the interaction between a regulator and multiple service providers: the regulator provides path-based subsidies at the upper level, while non-cooperative service providers set link prices to maximize their own profits at the lower level. \cite{pinto2020joint} adopted a bilevel structure to address a resource allocation problem, designing a multimodal system with autonomous vehicle (AV) fleets while capturing time-dependent mode choice behavior. Similarly, \cite{bandiera2024mobility} formulated a bilevel problem in which service operators determine optimal strategies at the upper level and traveler behaviors are modeled at the lower level. Instead of a centralized operator, each service provider maximizes its own profit, resulting in a Nash equilibrium at the upper level.

The deterministic assignment game in \citep{liu2024demand, yao2024design} was extended to a probabilistic approach by \cite{liu2024modeling} to reflect stochastic coalition choice. The assignment model adopts a similar path-based route choice framework as the capacitated SUE model from \cite{bell1995stochastic}, but the flows reflect the joint choices of the travelers and operators on which coalitions to match. This difference from a conventional route choice model is reflected in the added operator utility term in the lower-level objective function. Because the conventional SUE model requires path-based solutions, path set generation is also a prerequisite of the solution process for the stochastic assignment game in \cite{liu2024modeling}. However, high quality path set generation is a non-trival task as discussed in \citep{prato2009route}. Furthermore, the logit-based SUE in the lower level of a bilevel optimization makes it hard to integrate into a single level problem for solving to global optimality. As such, the path-based stochastic assignment game is not well-suited for large-scale applications.

\subsection{Mobility Hubs}
Mobility hubs can serve as key gateways that facilitate multi-modal mobility services. However, studies related to mobility hubs mainly focus on facility location and resource allocation, i.e. aspects of traditional transit hub design in which interactions between operators are ignored. Previous research has primarily focused on locating mobility hubs to improve accessibility and reduce social inequity, with the involved facilities typically centered around a single mode. \cite{caggiani2020approach,caggiani2020equality} proposed models to optimally place bike-share stations to enhance spatial fairness in mobility access. Similarly, \cite{duran2021demand} developed a model to site bike-share stations with the goal of minimizing spatial inequity. Other studies on hub location decisions also treat social welfare as a central research focus \citep{banerjee2020optimal, aydin2022planning}. \cite{frank2021improving} improve rural accessibility by using mobility hubs as gateways for multimodal trip itineraries. \cite{ren2025data} integrate synthetic data with observed data to evaluate the impacts of a MH on traveler disutilities.

Although numerous models have been proposed for designing mobility hubs, only a few studies consider the interactions among mobility service participants. \cite{nair2014equilibrium} was a pioneer in developing hub location models that incorporate traveler mode choices. They proposed a bilevel framework, with the hub operator as the upper-level decision maker and travelers as the lower-level decision makers. \cite{ma2019dynamic} provide a ridesharing operation strategy that prioritize interconnections with transit networks. Though not specifically designed for mobility hub application, the case study uses LIRR stations as key points for deploying proposed ridepooling dispatch and fleet repositioning strategy. The results show that the presence of a trunk commuter transit network would redistribute optimal ridepooling services to the last mile to increase effective service capacity by mitigating the "Wild Goose Chase" phenomenon (see also \cite{ouyang2023measurement}). \cite{xanthopoulos2024optimization} also demonstrate the importance of considering user preferences in determining optimal hub locations and capacities. Their framework decomposes hub decisions and traveler choices into multiple modules, and a customized metaheuristic was developed to apply the framework at a city scale. 

\subsection{Research gaps and our contributions}
Although numerous studies have developed models for designing multimodal service networks and mobility hubs, explicitly integrating mobility hubs into MaaS optimization frameworks remains a significant and under-explored area. \cite{nair2014equilibrium} and \cite{xanthopoulos2024optimization} captured the interaction between operator resource allocation decisions and traveler path choices. However, such analyses primarily treat hubs as mere transfer points, without explicitly considering their strategic role in operator decision-making processes. Specifically, service pricing and subsidy between operators are not considered key decision variables, even though they play crucial roles in influencing traveler choices. \cite{xi2024single, xi2024strategizing} proposed frameworks that leverage service capacity and pricing strategies in designing MaaS ecosystems. However, these frameworks do not consider the spatial aspects of resource allocation decisions and traveler itineraries, nor do they incorporate physical infrastructure when deciding service strategies. 

The frameworks proposed by \cite{liu2024demand} and \cite{yao2024design} effectively capture the complex interactions among the MaaS platform, service providers, and travelers while considering flexible service strategies in multimodal networks. Nonetheless, the role of mobility hubs as key levers in MaaS ecosystem design remains unaddressed. Additionally, the above-mentioned studies often require substantial computational resources, hindering their deployments on a city-level scale. These gaps call for further modeling approaches capable of capturing strategic interactions among operators and travelers, enhanced by more scalable computational methods. A path-based stochastic assignment game approach \citep{liu2024modeling} addresses some of the scalability and capturing of heterogeneous preferences, but faces others with path generation requirements.

This study addresses these gaps by rigorously integrating MH characteristics as platform control levers into MaaS market models, thereby contributing to both theoretical advancement and practical implementation strategies. The main contributions of this study are summarized as follows:

\begin{itemize}
    \item We propose a bilevel mathematical model for the general MaaS platform design problem, given the network structure. The model aims to maximize MaaS platform revenue (if private) or social welfare (if public) while capturing the two-sided matching between mobility hub costs and traveler utility at the lower level. The upper-level decision variables are the OD-based traveler-facing pricing that grants access to the MaaS services and operator-dedicated subsidies, while the lower-level decision variables are link flows and mobility service node capacities. The proposed model involves heterogeneous operator participation and the competition mechanism, using subsidies as a lever to transfer utility between travelers and operators with heterogeneous service elements.
    \item We adopt a random utility model (RUM) based formulation to capture the choices of the coalitions between traveler and operator choices in the lower level as proposed by \cite{liu2024modeling}. We extend it to a link-based formulation using the PURC approach proposed by \cite{fosgerau2022perturbed}. The approach offers improved scalability while still capturing link flows jointly determined by service operators and travelers, similar to other RUM formulations (e.g., stochastic user equilibrium, SUE).
    \item We reformulate the bilevel problem into a single-level problem using the Karush–Kuhn–Tucker (KKT) conditions of the lower-level problem. A gap function-based approach is proposed to accelerate the solution process. Because of the simplicity of the PURC-based quadratic program in the lower level, a global optimum can be attained within a reasonable time frame when applied to large-scale cases.
    % \item We extend the base model to include facility opening variables. Due to the increased complexity resulting from the introduction of binary variables for facility location decisions, we develop a customized matheuristic algorithm to efficiently obtain high-quality solutions.
    \item We apply the model to a real-world case study based on several LIRR stations and their surrounding areas to demonstrate the computational efficiency. Various operational schemes are also tested to demonstrate the impact of decisions such as pricing and subsidies. The case study provides insights for real-world policy-making in establishing mobility hubs to promote multimodal trips.
\end{itemize}

\section{Proposed mobility hub platform design methodology}

We first present the multimodal network structure and service assumptions. We then present the model formulations in detail. The model formulation for MH platform design takes the general stochastic MaaS assignment game from \citep{liu2024modeling} and re-formulates it within an equivalent PURC framework with explicit MH decision variables. The goal is to provide a scalable framework that jointly decides the service decisions from the operator side and path choices from the customer side, while the platform sets service subsidies. 

\subsection{Network structure and model assumptions}

The general notation for MaaS platform design used for the MH platform design follows a multicommodity network design (MCND) problem. The multimodal network $G$, denoted by $(N, A)$, consists of all nodes $N$ and links $A$, serving all origin-destination (OD) and segment groups $S$. Origins $O$ and destinations $D$ are aggregated into centroids $ \{O \cup D\} \subseteq N_S \subset N$. The services include a single fixed route transit (FT) operator acting as the platform for designing the MHs, and a set of MOD services that can choose to participate as feeders for these MHs through the platform. The model can be trivially modified to include fixed route buses run by other operators as feeders as well, but for simplicity in notation we assume this MH system deals only with a single fixed route platform/operator and one or more MOD feeder operators. The FT operator $f$ belongs to a singleton set $F$ and provides service on its own subnetwork $G_f$, which consists of node set $N_f \subseteq N_F \subset N$ and link set $A_f \subseteq A_F \subset A$. Similarly, each MOD fleet $m$ belongs to the MOD operator set $M$, with each operating a subnetwork $G_m$ consisting of node set $N_m \subseteq N_M \subset N$ and service link set $A_m \subseteq A_M \subset A$. A MOD service link represents a pickup at origin and drop-off at a MH, where the link attributes reflect average performance of the fleet serving this request under macroscopic modeling assumptions. Congestion effects, such as waiting for service, are captured by a node capacity at the origin. These "in-platform" links are store-and-forward links with service queues. Each MOD operator $m$ determines the capacity for nodes $i\in N_m$, which are modeled as capacitated virtual links $A_i^{\sim}$ to access the MOD subnetwork. The capacities are used to capture queue delays that arise from congestion while waiting for service. By representing node capacities as virtual links, we maintain a link-based network structure, ensuring consistency when estimating coefficients as per \cite{fosgerau2022perturbed}.

\begin{figure}
    \centering
    \includegraphics[width=0.75\linewidth]{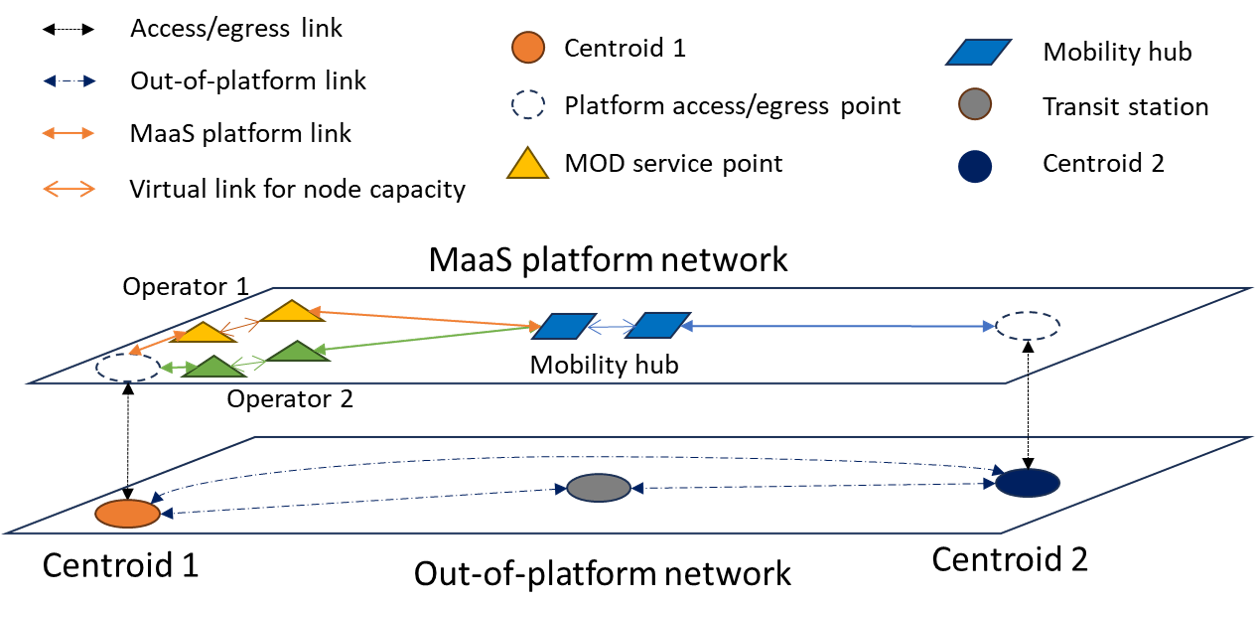}
    \caption{Network illustration for one OD pair}
    \label{fig:Network_illustration}
\end{figure}

A dummy subnetwork is used to represent travel outside of the MH platform's multimodal services (e.g., private car, telecommuting, or an alternative platform), i.e. "out-of-platform" links. For all traveler groups $S$, a set of platform access and egress links $A_S^+$ and $A_S^-$ are used to connect respective origins and destinations with the MH based platform. \textbf{Fig.~\ref{fig:Network_illustration}} illustrates the multimodal network structure for one OD pair. Travelers choose their preferred paths encompassing combinations of services across the network for each OD pair segment $s \in S$, which may use the platform to access the hub and taken the fixed route transit to the destination (and vice versa), use their own modes to get to the transit station to go to the destination, or use their own mode entirely to get to the destination. 

In the case of a mobility hub-oriented system, we consider a FT operator also serving as a platform regulator, working with MOD operators to provide services at designated MHs through a subset of nodes $H \subset N$. These operators may be micromobility providers, ridehail, or microtransit services. The hubs in this setting represent cyberphysical gateways in which participating MOD operators that pick up or drop off passengers purchasing trip bundles within these geofenced hubs could be subsidized. We introduce MOD feeder services centered at each hub $h \in H$, with direct access to surrounding nodes $i \in N_M$ represented by links $A_{i,h}^+$. The platform operator regulates total inbound flow to ensure it does not exceed hub capacity. For consistency, these hub capacities are modeled as virtual links $A_h^{\sim}$.

For travelers using the hub centered platform services, an OD-based access fee is charged at the access link $A_{s}^{+}$ for $s \in S$. After paying the access fee, traversing links belonging to the platform does not incur additional ticket charges as long as the traveler does not exit the platform-based service during the trip.

As shown in \textbf{Fig.~\ref{fig:Network_illustration}}, we separate subnetworks to be "in-platform" and "out-of-platform". We focus on the OD-based pricing at the access link for travelers and service strategy decisions of operators in-platform, including operator dedicated subsidies and service capacities. Since out-of-platform includes all other options out there (including driving, walking, telecommuting, or using MOD services via pay-as-you-go option), they are assumed to be uncapacitated and we do not explicitly model competition between the designed platform and other external modes. Instead, we focus on the cooperative game between operators in-platform.

The in-platform services can be provided by multiple operators. Travelers can choose their preferred service links after entering the in-platform subnetwork by paying the access price, which is a decision variable in the model. While inside the subnetwork, only non-fee related costs are induced. For participating operators, their operating costs are covered by the subsidy provided by the platform. Other decision variables include the service capacities at either the MH or the MOD service nodes within the in-platform subnetworks. Table \ref{tab:variables} provides the full list of notations. In addition, we list the model assumptions as follows. To avoid confusion, we define travel cost as the non-monetary penalties associated with links such as travel time. Access price is the monetary values required for entering the in-platform subnetwork. Both terms are converted to utilities.

\begin{table}[htbp]
\centering
\caption{Model Variables and Parameters}
\begin{tabular}{ll}
\hline
\textbf{Notation} & \textbf{Description} \\
\hline
\multicolumn{2}{l}{\textit{Sets}} \\
\hline
$N$ & Set of all nodes in the network \\
$A$ & Set of all links in the network \\
$S$ & Set of origin-destination (OD) pairs and population segments \\
$N_S \subset N$ & Set of centroid nodes (OD origins/destinations) \\
$F$ & Set of fixed route transit (FT) operators \\
$N_F \subseteq N$ & Node set of all FT operators \\
$A_F \subseteq A$ & Link set of all FT operators \\
$N_f \subseteq N_F$ & Node set for FT operator $f \in F$ \\
$A_f \subseteq A_F$ & Link set for FT operator $f \in F$ \\
$M$ & Set of MOD (Mobility on Demand) operators \\
$N_M \subseteq N$ & Node set of all MOD operators \\
$A_M \subseteq A$ & Link set of all MOD operators \\
$N_m \subseteq N_M$ & Node set for MOD operator $m \in M$ \\
$A_m \subseteq A_M$ & Link set for MOD operator $m \in M$ \\
% $N_O$ & Node set for out-of-platform travel (e.g., private car) \\
$A_O$ & Link set for out-of-platform travel \\
$A_0$ & Set of transfer links between centroids and service nodes, and between layers \\
$H \subset N$ & Set of tentative mobility hubs \\
$A_H$ & Set of added links connecting hubs and service nodes \\
$A_h^{\sim}$ & Set of virtual hub links $h \in H$ with capacity\\
$A^-_i$ & Set of outbound links of node i \\
$A^+_i$ & Set of inbound links of node i \\
$A_{i}^{\sim}$ & Set of virtual service node links $i\in N_M$ with capacity\\
$A_{i,h}^+$ & Set of direct links from MOD node $i\in N_M$ to hub node $h$ \\
$A_S^+$ & Set of access links entering MH based in-platform subnetwork \\

\hline
\multicolumn{2}{l}{\textit{Input Variables (Parameters)}} \\
\hline
$d_l$ & Length of link $l \in A$ \\
$c^{t}_{l}$ & Traveler cost per unit of link $l \in A$ \\
$c^{o}_l$ & Operator cost per unit of link $l \in A$ \\
$q_{s}$ & OD demand for each OD pair $s \in S$ \\
$z_l$ & Capacity of access links towards MOD nodes $l \in A_i^{\sim}, i\in N_M$ \\
$c_i$ & Per capacity cost of MOD node $i \in N_M$ \\
$a_{i,l}$ & Node-link incidence: $-1$ if node $i$ is origin of link $l$, $+1$ if destination, $0$ otherwise \\
$\hat{p}_l$ & Cap of access price for link $l \in A_S^+$ \\
$\bar{q}$ & Average OD demand \\
$\alpha_1$ & Weight of traveler utility\\
$\alpha_2$ & Weight of operator utility\\

\hline
\multicolumn{2}{l}{\textit{Decision Variables}} \\
\hline
$x_{l,s}$ & Normalized flow on link $l \in A$ for OD pair $s \in S$ \\
$v_l$ & Proportion of maximum capacity opened for access link $l \in A_{0i}$ \\
% $y_i$ & Opening decision of MOD node $i \in N_M$, binary \\
$p_{l,s}$ & Service price on platform access link $l \in A_s^{+}$ for each OD pair $s \in S$ \\
$r_l$ & Operator subsidy provided by platform based on link $l \in A_i^{\sim}$ \\
$b_l$ & Capacity of hub virtual link $l\in A_h^{\sim}$ \\
\hline
\end{tabular}
\label{tab:variables}
\end{table}

List of assumptions:
\vspace{-2ex}
\begin{itemize}
    \item The proposed model only considers strategic planning perspective. Therefore, the model only involves static flow assignment.
    \item The MH platforms are centrally controlled, which involves the subsidy decisions.
    \item All cost-related terms are link-additive.
    \item Utilities are transferable between travelers, operators, and the platform.
    \item Travel cost of the in-platform links, FT service links, and out-of-platform links are fixed.
    \item Transfer links between centroids and FT nodes have fixed travel cost including wait and walk time, which represent the overall cost of accessing FT services.
    \item MOD access links are store-and-forward links with service queues, meaning that the access flow remains uncongested until capacity is reached. Once capacity is exceeded, excess flows are diverted to other links with available capacity. The queue delay is quantified by the capacity Lagrange multiplier.
    \item All FT service links and links in the out-of-service subnetwork are uncapacitated. Congestion from background traffic are directly incorporated into the travel costs. If an out-of-service link reflects multiple modes, its travel cost reflects the expected minimized disutility option of those modes. Alternatively, each parallel mode may be modeled separately.
    \item Operating costs of MOD services consist of service link costs and access capacity allocation costs.
    \item Operator costs for service links, MOD nodes, and mobility hubs scale linearly with demand. We acknowledge that nonlinear functions, such as the Cobb–Douglas form \citep{zhang2020modeling}, are often used to capture operation cost in mobility services. However, given that this model addresses macroscopic system planning under steady-state conditions, the linear assumption provides a necessary balance between model accuracy and computational efficiency.
\end{itemize}

\subsection{Lower-level traveler-operator coalition choice assignment model using PURC based structure}
The MH platform design problem is adapted from the stochastic MaaS assignment game model in \cite{liu2024modeling}. The assignment game is decomposed into a bilevel problem where the upper level involves platform decisions, such as pricing and subsidies, while the lower level involves the joint decisions of the travelers and mobility operators. In the adaptation to the MH platform design problem, the platform makes decisions about how much price to charge travelers for each OD segment and how much to subsidize the MOD operators in the upper level, while the lower level determines the MH and MOD capacities and traveler flows.

Instead of using the path-based SUE-style formulation in \cite{liu2024modeling}, we adapt the link-based PURC model proposed by \cite{fosgerau2022perturbed} to capture the probabilistic matching mechanism. Eq.~(\ref{eq:service_model}) is the link-based stochastic assignment problem $L_1$ based on the PURC concept for a MH platform.

\vspace{-2ex}
\begin{subequations}\label{eq:service_model}
\begin{align}
L_1\min_{\{x_{l,s},\, v_l\}} \quad
    & \Phi_1 = \sum_{l \in A} \sum_{s \in S} \bar{q} d_l x_{l,s}^2 \notag \\
    & \quad + \alpha_1\bar{q}\sum_{s \in S}\left(\sum_{l \in A_s^{+}} d_l p_{l,s} x_{l,s} + \sum_{l \in A/A_{s}^{+}}d_lc_lx_{l,s}\right) \notag \\
    & \quad + \alpha_2\left(
        \sum_{l \in A} \sum_{s \in S} d_l c^{o}_l x_{l,s} q_s
        +\sum_{l \in A_i^{\sim}, i\in N_M} \left(z_l c_i v_l - \sum_sr_lq_sx_{l,s}\right)
    \right) \label{eq:service_model_obj} \\
\text{s.t.}\quad
    & \sum_{l \in A} a_{i,l} x_{l,s} =
    \begin{cases}
    -1, & \forall i = o , \\
    1, & \forall i = d, \\
    0, & \forall i \in N \setminus \{o,d\}
    \end{cases} \forall s = (o,d) \in S,\, \forall i \in N
    \label{eq:service_model_flow}\\
    & \sum_{s \in S} x_{l,s} q_{s} \leq z_l v_l
    \quad \forall l \in A_i^{\sim}, \forall i \in N_M
    \label{eq:service_model_modcap}\\
    & \sum_{l \in A^+_{i, h}} \sum_{s \in S} x_{l,s} q_{s} \leq b_{\tilde{l}}
    \quad \forall \tilde{l} \in A_h^{\sim}, \forall h \in H
    \label{eq:mh_modcap}\\
    & 0 \leq  x_{l,s} \leq 1
    \quad \forall l \in A,\, \forall s \in S
    \label{eq:service_model_xnonneg}\\
    & 0 \leq v_l \leq 1
    \quad \forall l \in A_i^{\sim}, \forall i \in N_M
    \label{eq:service_model_vnonneg}
\end{align}
\end{subequations}

The assignment model $\Phi_1$ serves as a lower level problem to the corresponding assignment game, but can also be run as an independent model for determining a multimodal network design problem with stochastic assignment, such as determining the welfare impact of an existing pricing design of a traditional transit hub without MH subsidy. 

The decision variables are the percent link flow assignment $x_{l,s}$ and percent of maximum MOD service capacity $v_l$ assigned to the access links. The platform access price $p_{l,s}$ is the decision variable in the upper level model, which is presented in the following subsection. In $\Phi_1$, $p_{l,s}$ is treated as an input. Eq.~(\ref{eq:service_model_flow}) ensures the link flow conservation. For each OD pair $s=(o,d)$ with $o$ being the origin and $d$ being the destination, the source node $o$ initiates all normalized flow, and the sink node $d$ terminates all normalized flow. Eq.~(\ref{eq:service_model_modcap}) controls the capacity and congestion effect of the platform service node. All access links flowing into the service nodes are bounded by the assigned capacity the operator determines. Eq.~(\ref{eq:mh_modcap}) is dedicated to mobility hub facilities. For all service flows ending at mobility hubs, the sum of flow needs to be capped by the assigned mobility hub capacity, which is determined by the platform in the upper level model. Eqs.~(\ref{eq:service_model_xnonneg} - \ref{eq:service_model_vnonneg}) are the bounds for the two sets of decision variables. 

The objective consists of three terms: the perturbed utility, traveler utility, and the operator utility. Introduced by \cite{fosgerau2022perturbed}, the PURC model follows the general RUM structure with link-based property to capture the dispersion of flows across a network. The flow assignment variable $x_{l,s}$ represents the probability of choosing link $l$ for travelers in OD pair $s$. The first term in Eq.~(\ref{eq:service_model_obj}) depicts the dispersion effect. Without this term, the model becomes an all-or-nothing assignment problem. By adding a strictly increasing function where the first and second derivatives becomes zero at the origin, a dispersion effect is added that permits deviation in route choice. \cite{fosgerau2022perturbed} demonstrate that both quadratic and entropic functions are effective. We select the quadratic form for this study due to its simplicity, and it maintains a high level of accuracy as shown in \cite{fosgerau2022perturbed}. The second term is the sum of traveler utility, which depicts the sum of all link utilities. We simplify the term by only involving the platform access price element and general travel cost element. It can be expanded to involve all utility related elements such as externalities, as long as they can be converted to utility values. For a comprehensive explanation of the disperse effect and the utility features that can be included, readers can refer to \cite{fosgerau2022perturbed} for more information. The original PURC objective is constructed by combining the first and second terms. We show this in \textbf{Proposition}~\ref{prop_1}.

\begin{proposition}\label{prop_1}
The first and second terms combined in Eq.~(\ref{eq:service_model_obj}) are equivalent to the PURC objective proposed by \cite{fosgerau2022perturbed}.
\end{proposition}

\begin{proof}
The first and second term of Eq.~(\ref{eq:service_model_obj}) can be written as Eq.~(\ref{eq:rewritten_obj}).
\begin{equation} \label{eq:rewritten_obj}
    \min_{\{x_{l,s},\, v_l,\}} \quad \Phi_1 =\bar{q}\left(\sum_{l \in A} \sum_{s \in S} d_l x_{l,s}^2  + \alpha_1\sum_{s \in S}\left(\sum_{l \in A_s^{+}} d_l p_{l,s} x_{l,s} + \sum_{l \in A/A_{s}^{+}}d_lc_lx_{l,s}\right)\right)
\end{equation}
The coefficient $\bar{q}$ is a constant. The quadratic term is the perturbed term $F(x)$, and the linear term is the link-additive utility term $U(x)$. Both terms are weighted by traversed link length $d_l$. The structure is therefore identical to the PURC objective proposed in \cite{fosgerau2022perturbed} by multiplying with a constant $\bar{q}$.
\end{proof}

We introduce the third term to form the coalition between travelers and operators in a similar fashion to the model proposed by \cite{liu2024modeling}. In \textbf{Proposition}~\ref{prop_2}, we show that $\Phi_1$ provides a flow assignment model with a traveler-operator coalition under the rule of PURC.

\begin{proposition}\label{prop_2}
$\Phi_1$ yields the link flow assignment under PURC framework with traveler and operator coalition choice.
\end{proposition}

\begin{proof}
We divide Eq.~(\ref{eq:service_model_obj}) by the constant term $\bar{q}$ and write it as Eq.~(\ref{eq:lower_obj_rewritten}).
\begin{equation} \label{eq:lower_obj_rewritten}
\begin{aligned}
    \min_{\{x_{l,s},\, v_l\}} \quad
    \Phi'_1 ={}& \sum_{l \in A} \sum_{s \in S} d_l x_{l,s}^2  \\
    &+ \alpha_1 \sum_{s \in S} \left( \sum_{l \in A_s^{+}} d_l p_{l,s} x_{l,s}
           + \sum_{l \in A\setminus A_{s}^{+}} d_l c_l x_{l,s} \right) \\
    &+ \alpha_2 \left(
        \sum_{l \in A} \sum_{s \in S} d_l c^{o}_l x_{l,s} \frac{q_s}{\bar{q}}
        + \sum_{l \in A_i^{\sim}} \frac{ z_l c_i v_l - \sum_{s} r_l q_s x_{l,s} }{\bar{q}}
    \right)
\end{aligned}
\end{equation}

The Lagrangian of $\Phi'_1$ is written as Eq.~(\ref{eq:lagrangian_phi_1}).
\begin{equation} \label{eq:lagrangian_phi_1}
\begin{split}
    \min_{\{x_{l,s},\, v_i,\}} \quad
    L'_1 ={}& \sum_{l \in A} \sum_{s \in S} d_l x_{l,s}^2  \\
    &+ \alpha_1 \sum_{s \in S} \left( \sum_{l \in A_s^{+}} d_l p_{l,s} x_{l,s}
           + \sum_{l \in A\setminus A_{s}^{+}} d_l c_l x_{l,s} \right) \\
    &+ \alpha_2 \left(
        \sum_{l \in A} \sum_{s \in S} d_l c^{o}_l x_{l,s} \frac{q_s}{\bar{q}}
        + \sum_{l \in A_i^{\sim}} \frac{ z_l c_i v_l - \sum_{s} r_l q_s x_{l,s} }{\bar{q}}
    \right) \\
    & + \sum_{s\in S}\sum_{i\in N}\mu_{i,s}(\sum_{l}a_{i,l}x_{l,s}-f_{i,s})\\
    & + \sum_{l\in A_i^{\sim}, i\in N_M}\lambda_l(\sum_{s \in S} x_{l,s} q_{s} - z_l v_l)\\
    & + \sum_{\tilde{l} \in A_h^{\sim}, h\in H}\varphi_{\tilde{l}}(\sum_{l\in A_{i,h}^+}\sum_{s\in S}x_{l,s}q_s - b_{\tilde{l}}) \\
    & + \sum_{l\in A} \sum_{s \in S} \beta_{l,s}(x_{l,s}-1)+\sum_{l\in A_i^{\sim}, i\in N_M}\pi_{l}(v_l-1)
\end{split}
\end{equation}

The optimal solution of $x_{l,s}$ is obtained by having $\frac{\partial L'_1}{x_{l,s}}=0$ for all chosen links $l$ for OD pair $s$.

\begin{equation} \label{eq:derive_1}
\begin{split}
    \frac{\partial L'_1}{x_{l,s}} &=2d_l x_{l,s}+ d_l \left(\alpha_1\left(p_l +  c_l\right) + \alpha_2 \frac{(c^{o}_l-r_l) q_s}{\bar{q}}\right)\\
    & + \sum_{i\in {l^+,l^-}}\mu_{i,s}a_{i,l} +\lambda_l q_{s} + \sum_{\tilde{l}^-\in l^+}\varphi_{\tilde{l}} q_s +\beta_{l,s} = 0
\end{split}
\end{equation}

We slightly change the notation of the term associated with $\alpha_1$ in Eq.~\ref{eq:derive_1} for simplicity. By setting $c_l$ to zero when $l\in A_s^+$ and $p_l$ to zero when $l\in A\setminus A_s^+$, Eq.~\ref{eq:derive_1} is equivalent to the first order derivative of Eq~\ref{eq:lagrangian_phi_1}. Because the chosen links will have strictly positive flows, the $v_l$ associated with those links also fulfill the optimality conditions.

\begin{equation} \label{eq:derive_2}
    \frac{\partial L'_1}{v_{l}} = \alpha_2 \frac{z_l c_i}{\bar{q}} - \lambda_{l} z_l + \pi_l = 0
\end{equation}

By re-writing Eq.~(\ref{eq:derive_2}), we obtain $\lambda_l= \frac{\alpha_2 c_i}{\bar{q}}+\frac{\pi_l}{z_l}$. By plugging it in Eq.~(\ref{eq:derive_1}), we have Eq.~(\ref{eq:derive_3}).

\begin{subequations} \label{eq:derive_3}
\begin{align}
    \frac{\partial L'_1}{x_{l,s}} &=2d_l x_{l,s}+ d_l \left(\alpha_1\left(p_l +  c_l\right) + \alpha_2 \frac{(c^{o}_l-r_l) q_s}{\bar{q}}\right) \\
    & + \sum_{i\in {l^+,l^-}}\mu_{i,s}a_{i,l} +\sum_{i\in l^+}\alpha_2 c_i\frac{q_{s}}{\bar{q}} + \frac{\pi_lq_s}{z_l} + \sum_{l'^-\in l^+}\varphi_{l'} q_s + \beta_{l,s} = 0
\end{align}
\end{subequations}

Eq.~(\ref{eq:derive_3}) can be rewritten as Eq.~(\ref{eq:derive_4}).

\begin{equation} \label{eq:derive_4}
   x_{l,s} = -  \frac{1}{2}\left(\alpha_1\left(p_l +  c_l\right) + \alpha_2 \frac{(c^{o}_l-r_l) q_s}{\bar{q}}\right) - \frac{1}{2d}(\sum_{i\in {l^+,l^-}}\mu_{i,s}a_{i,l} +\sum_{i\in l^+}\alpha_2 c_i\frac{q_{s}}{\bar{q}} +\frac{\pi_lq_s}{z_l} + \sum_{l'^-\in l^+}\varphi_{l'} q_s + \beta_{l,s})
\end{equation}

Eq.~(\ref{eq:derive_4}) depicts the flow assignment decisions for all links $l$ with positive flows for each OD pair $s$. The first term consists of the link disutility from both travelers and operators. The link disutilities from the two sides are weighted by $\alpha_1, \alpha_2$, and $q_{s}/\bar{q}$. $\alpha_1$ and $\alpha_2$ are traveler and operator dedicated coefficients that provide different weights towards the overall utility objective. The ratio of $\alpha_1/\alpha_2$ depicts which side having higher impact towards the flow assignment behavior. $q_{s}/\bar{q}$ is an additional weight that considers the impact from different OD pair volumes. The flow assignment variable $x_{l,s}$ is considered as the possibility of choosing link $l$ for OD group $s$. Each OD pair shall be treated equally from a traveler perspective to reflect stochastic user-equilibrium (SUE) behavior from the perspective of the coalitions of travelers and operators. The weight $q_{s}/\bar{q}$ is the lever that places more focus on OD pairs with higher OD volume. 

The second term consists of node and capacity specific values. The first element represents the magnitude of attraction of link $l$ in OD pair $s$, which is measured by the difference between the source and sink node Lagrange multipliers. The second element represents the disutility from the cost of serving a node for MOD operators. For links that do not belong to the access link set leading to MOD nodes $N_M$, the second element is dropped. For the the access links leading to MOD nodes, the node cost disutility from the operator side is counted towards the flow assignment, which is also weighted by $q_{s}/\bar{q}$. The remaining three terms are the Lagrange multipliers associated with MOD service node capacity, hub capacity, and link flow constraints.

Eq.~(\ref{eq:derive_4}) mimics the structure of the optimal link-flow decision variable with non-zero value shown in \cite{fosgerau2022perturbed}. This concludes the proof.
\end{proof}

Representing the differences in the joint behavior of the traveler and operator coalitions, $\alpha_1$ and $\alpha_2$ would typically be fitted to data in practice to help regulate the dispersion effect observed in the data. As shown in \cite{fosgerau2022perturbed}, route dispersion intensifies when the first term in Eq.~(\ref{eq:service_model_obj}) holds greater relative weight. If the observed flows have less divergent route choice behavior, the magnitudes of $\alpha_1$ and $\alpha_2$ can be both raised while maintaining the same ratio of $\alpha_1/\alpha_2$.

\cite{bell1995stochastic} proved that the congestion effect can be reflected by the capacity Lagrange multiplier in a capacitated network, which is also shown in \cite{liu2024modeling}. We follow the same logic to show that the congestion effect can also be evaluated by capacity Lagrange multipliers in the PURC based assignment model.

\begin{proposition}\label{prop_3}
The congestion effect of waiting for service at a capacitated MOD access link $l\in A_i^{\sim}, i \in N_M$ is captured by the corresponding Lagrange multiplier $\lambda_{l}$. Similarly, the congestion effect caused by the capacitated hub facility $l\in A_h^{\sim}, h \in H$ is captured by the corresponding Lagrange multiplier $\varphi_{l}$.
\end{proposition}

\begin{proof}
We use the congestion effect at MOD access links as an example. Assume an access link $l$ with sink node $i$ is at capacity ($v_l = 1$). If there is an additional flow $\delta_s$ for OD pair $s$ on link $l$, the capacity constraint Eq.~(\ref{eq:service_model_modcap}) is therefore exceeded by $\frac{\delta_{s}q_s}{z_l}$. Assume the additional flow switches to link $l'$ while the chosen downstream links have identical costs compared to the chosen downstream links sourcing from node $i$. The change of objective value measured by Eq.~(\ref{eq:service_model_obj}) is written as Eq.~(\ref{eq:congestion}).

\begin{equation} \label{eq:congestion}
   d \Phi_1 = \bar{q} (d_{l} (2x_{l,s}+\delta_s) -d_{l'} (2x_{l',s}+\delta_s)) + \bar{q} \delta_s (d_l p_l - d_{l'}p_{l'}+ d_l c^{t}_l -d_{l'} c_{l'}^t))
    + \alpha \delta_s(
         d_l c^{o}_l - d_{l'} c^{o}_{l'})q_s
\end{equation}

As shown in Eq.~(\ref{eq:congestion_unit}), when the right-hand side of the node capacity constraint Eq.~(\ref{eq:service_model_vnonneg}) increases by 1, the objective value of the KKT in Eq.~(\ref{eq:lagrangian_phi_1}) increases by $\pi_l$. Therefore, we can rewrite Eq.~(\ref{eq:congestion}) as Eq.~(\ref{eq:congestion_unit_rewrite}).

\begin{equation} \label{eq:congestion_unit}
   \frac{d\Phi_1}{\delta_{s}q_s/z_l} = \frac{z_l \bar{q}}{q_s} (d_{l} (2x_{l,s}+\delta_s) -d_{l'} (2x_{l',s}+\delta_s)) +\frac{z_l \bar{q}}{q_s} (d_l p_l - d_{l'}p_{l'}+ d_l c^{t}_l -d_{l'} c_{l'}^t))
    + \alpha z_l(
         d_l c^{o}_l - d_{l'} c^{o}_{l'})
         = \pi_l
\end{equation}

\begin{equation} \label{eq:congestion_unit_rewrite}
  \frac{\bar{q}}{q_s} (d_{l} (2x_{l,s}+\delta_s) -d_{l'} (2x_{l',s}+\delta_s)) +\frac{\bar{q}}{q_s} (d_l p_l - d_{l'}p_{l'}+ d_l c^{t}_l -d_{l'} c_{l'}^t))
    + \alpha (
         d_l c^{o}_l - d_{l'} c^{o}_{l'})=\frac{\pi_l}{z_l}
\end{equation}

If link $l$ and $l'$ both having 1 unit of length, Eq.~(\ref{eq:congestion_unit_rewrite}) is equivalent to Eq.~(\ref{eq:congestion_effect}).

\begin{equation} \label{eq:congestion_effect}
   x_{l,s}-x_{l',s}=-(p_l - p_{l'})- (c^{t}_l + \frac{\pi_l q_s}{z_l\bar{q}}-c_{l'}^t)
    - \alpha (
         c^{o}_l - c^{o}_{l'})\frac{q_s}{\bar{q}}
\end{equation}

Eq.~(\ref{eq:congestion_effect}) shows that the congestion effect caused by the capacitated node can be evaluated by the Lagrange multiplier of the node capacity constraint. When the switched links are identical in length, the congestion effect equals to $\frac{\pi_l q_l}{z_l\bar{q}}$. The same logic can be applied to capacitated hub links. This concludes the proof.
\end{proof}

The lower-level assignment model is a convex quadratic programming (QP) problem, which can be easily solved with off-the-shelf solvers at scale. Moreover, it exhibits characteristics that allow us to solve the bilevel problem to global optimality as well.

\subsection{Upper level model formulation}

We present the upper-level problem in this subsection. The decision-maker of the upper-level problem is the MH platform, which determines the optimal OD-based pricing strategy for accessing the mobility service, the capacity assigned to hub locations, and the subsidy provided to incentivize operator participation. The bilevel structure follows the Stackelberg equilibrium principle, with the MaaS platform as the leader that sets pricing, capacity, and subsidy decisions, while travelers and operators are followers who make decisions to minimize their disutilities. We use revenue maximization as the upper-level objective. However, other equity or social welfare focuses objective can be applied for the upper level objective. From a revenue maximization perspective, the mobility hub design problem can be formulated as Eq.~(\ref{eq:upper_level}). 
\begin{subequations}\label{eq:upper_level}
\begin{align}
    \max_{\mathbf{p, r, b}} \quad &\Phi_{0} = \sum_{s}\sum_{l\in A_{s}^+}p_{s} \, x_{l, s} \, q_{s} - \sum_{l\in A_H}\sum_{i\in l^-}c_ib_l - \sum_{m\in M}\sum_{i \in N_m}\sum_{l\in A_i^{\sim}}\sum_{s\in S}r_l x_{l, s} q_s \label{eq:upper_obj} \\
    \text{s.t.} \quad & \mathbf{x} = \arg\min_{\mathbf{x}} \; \Phi_{1}(\mathbf{p, r, b}) \label{eq:lower_level} \\
    & 0 \leq p_{s} \leq \hat{p}, \quad \forall s\in S \label{eq:price_cap} \\
    & 0 \leq b_l \leq \hat{b}, \quad \forall l \in A_H \label{eq:cap_cap} \\
    & 0 \leq r_l \leq \hat{r}, \quad \forall l \in A_i^{\sim}, \forall i \in N_M \label{eq:sub_cap} \\
    & \sum_{i \in N_m}\sum_{l\in A_i^{\sim}}\sum_{s\in S}r_l x_{l, s} q_s - \sum_{l \in A_m}\sum_{s\in S}d_lc_l^o x_{l,s}q_s - \sum_{i\in N_m}\sum_{l\in A_i^{\sim}}z_lc_iv_l \geq 0, \quad \forall m\in M \label{eq:cost_floor}
\end{align}
\end{subequations}
The upper-level problem finds the optimal pricing strategy, along with the hub capacity and operator subsidy, that maximize the total revenue of the MH platform, while the service flows and MOD route capacities are determined by the lower-level problem. The first and second terms represent the total income from passenger access and the hub operating cost. The third term represents the total subsidy provided to participating MOD operators. The model assumes that $\hat{p}, \hat{b}$, and $\hat{r}$ reflect the existing conditions that the market accepts. These parameters could be further segmented for specific traveler and operator groups. We adopt uniform, system-wide values for simplicity in Eq.~(\ref{eq:price_cap} - \ref{eq:sub_cap}). Eq.~(\ref{eq:cost_floor}) requires that the subsidy for each MOD operator cover the operating cost so participation in the MH-oriented platform is profitable.

We can also slightly modify the objective so that subsidy also comes from outside resources (e.g., local authority). Assuming there is a fixed subsidy package $R$ provided to the platform operator, then by directly adding $R$ to Eq.~(\ref{eq:upper_obj}), we can effectively capture the impact of outside subsidy without modifying the constraints.

Aside from the PURC modifications, the bilevel problem presented in these two subsections feature several nontrivial differences from the generic MaaS assignment game model from \citep{liu2024demand} and \citep{liu2024modeling}. These include:
\vspace{-2ex}
\begin{itemize}
    \item In-platform links only include the ones inbound and outbound from MHs; the consequence is the capacity designs focus only on which last-mile service regions the MOD operators would serve (capacities approaching zero suggests not serving those regions) and the capacities for the MHs (capacities approaching zero here suggests closing the MH).
    \item There are out-of-platform links not just connecting centroid OD pairs, but also connecting centroids to the FT operators' stations as out-of-platform access links for travelers choosing to use the FT service without engaging with the platform.
    \item Whereas the stochastic assignment game in \citep{liu2024modeling} lacks a subsidy and the deterministic one in \citep{liu2024demand} only subsidizes a system that has no stable outcome otherwise, the proposed model determines both price and subsidy as upper-level decisions with a guaranteed stable outcome. As shown in \citet{liu2024modeling}, the stochastic assignment game is guaranteed to have a stable outcome with a single cost allocation corresponding to the optimal flow because of the endogeneity in the coalition choices. The added subsidy further assures the participating operator's profitability while strictly following the coalition structure, making the model more applicable in real-world adaptation.
    \item The general MaaS platform assignment game assumes there is no preexisting network, so as costs for operators go up, the optimal solution will have everybody move to out-of-platform alternatives. In the proposed model, there is a preexisting FT service, so if costs of operating MHs are too high, out-of-platform flows imply travelers either going to their destinations or to the transit stations on their own modes via driving, walking, or pay-as-you-go options. 
\end{itemize}

\subsection{Solution method}
The bilevel programming problem is nonconvex due to the dependency between leader and follower. The bilevel problem can be reformulated into a single-level constrained optimization problem by incorporating the KKT conditions of Eq.~(\ref{eq:service_model}) as additional constraints into Eq.~(\ref{eq:upper_level}). The reformulated single-level problem is written as Eq.~(\ref{eq:single_level}).
\begin{subequations}\label{eq:single_level}
\begin{align}
    \max_{\mathbf{p,x,v, \mu, \lambda, \beta, \pi}} \quad &\Phi'_{0} = \sum_{s}\sum_{l\in A_{s}^+}p_{s} \, x_{l, s} \, q_{s} - \sum_{l\in A_H}\sum_{i\in l^-}c_ib_l - \sum_{m\in M}\sum_{i \in N_m}\sum_{l\in A_i^{\sim}}\sum_{s\in S}r_l x_{l, s} q_s \label{eq:single_obj} \\
    \text{s.t.} \quad & 0 \leq p_{s} \leq \hat{p}, \quad \forall s\in S \label{eq:price_cap_single}\\
    & 0 \leq b_l \leq \hat{b}, \quad \forall l \in A_H \label{eq:cap_cap_single} \\
    & 0 \leq r_l \leq \hat{r}, \quad \forall l \in A_i^{\sim}, \forall i \in N_M \label{eq:sub_cap_single} \\
    & \sum_{i \in N_m}\sum_{l\in A_i^{\sim}}\sum_{s\in S}r_l x_{l, s} q_s - \sum_{l \in A_m}\sum_{s\in S}d_lc_l^o x_{l,s}q_s - \sum_{i\in N_m}\sum_{l\in A_i^{\sim}}z_lc_iv_l \geq 0, \quad \forall m\in M \label{eq:cost_floor_single}\\
    & x_{l,s}(2d_l x_{l,s}+ d_l \left(\alpha_1\left(p_l +  c_l\right) + \alpha_2 \frac{(c^{o}_l-r_l) q_s}{\bar{q}}\right) \nonumber \\
    & + \sum_{i\in {l^+,l^-}}\mu_{i,s}a_{i,l} +\sum_{i\in l^+}\alpha_2 c_i\frac{q_{s}}{\bar{q}} + \sum_{l'^-\in l^+}\varphi_{l'} q_s + \beta_{l,s}) = 0 \quad \forall l\in A, \forall s\in S \label{eq:link_kkt}\\
    & \alpha_2 \frac{z_l c_i}{\bar{q}} - \lambda_{l} z_l + \pi_l = 0, \quad \forall l \in A_i^{\sim}, \forall i\in N_M \label{eq:node_kkt} \\
    & \mu_{i,s}(\sum_l a_{i,l}x_{l,s} - f_{i,s})=0, \quad \forall i \in N, \forall s \in S \label{eq:node_slack}\\
    & \lambda_l(\sum_{s \in S} x_{l,s} q_{s} - z_l v_l)=0 \quad \forall l\in A_i^{\sim}, \forall i\in N_M \label{eq:cap_slack}\\
    & \varphi_{\tilde{l}}(\sum_{l\in A_{i,h}^+}\sum_{s\in S}x_{l,s}q_s - b_{\tilde{l}})=0, \quad \forall \tilde{l} \in A_h^{\sim}, \forall h \in H \\
    & \beta_{l,s}(x_{l,s} - 1) = 0, \quad \forall l \in A, s\in S \label{eq:link_slack} \\
    & \pi_l(v_l -1)=0, \quad \forall l\in A_i^{\sim}, \forall i\in N_M \label{eq:cap_slack_1}\\
    & \sum_{l \in A} a_{i,l} x_{l,s} =
    \begin{cases}
    -1, \quad \forall i = o , \\
    1, \quad \forall i = d, \\
    0, \quad \forall i \in N \setminus \{o,d\}
    \end{cases}
    \quad \forall s = (o,d) \in S,\, \forall i \in N
    \label{eq:single_level_flow}\\
    & \sum_{s \in S} x_{l,s} q_{s} \leq z_l v_l
    \quad \forall l \in A_i^{\sim}, \forall i \in N_M
    \label{eq:service_model_modcap_single}\\
    & \sum_{l \in A^+_{i, h}} \sum_{s \in S} x_{l,s} q_{s} \leq b_{\tilde{l}}
    \quad \forall \tilde{l} \in A_h^{\sim}, \forall h \in H
    \label{eq:mh_modcap_single}\\
    & \mathbf{p,x,v, \mu, \lambda, \beta, \pi} \geq 0 \label{eq:nonneg_single}
\end{align}
\end{subequations}

Eqs.~(\ref{eq:link_kkt}, \ref{eq:node_kkt}) are the first order conditions. Eq.~(\ref{eq:link_kkt}) is an augmented first order condition because of the zero-flow links involved in the final solution as indicated in \cite{fosgerau2022perturbed}. Eqs.~(\ref{eq:node_slack} - \ref{eq:cap_slack_1}) are the complementarity slackness conditions. It has been shown (e.g. \cite{dempe2002foundations}) that bilevel problems with linear upper level and convex quadratic lower level programs can be reformulated into a single level mathematical program with complementarity constraints (MPCC) that can be solved to global optimality. \cite{sinha2017review} provide an example algorithm based on branch-and-bound. 

However, the added KKT multipliers and complementarity constraints pose great scalability challenges, especially when facing a large network. Therefore, it is impractical to solve for a global optimal solution using branch-and-bound. To address this issue, we adopt the gap function-based method proposed by \cite{marcotte1996exact} to solve the bilevel problem effectively while still keeping high solution quality. Instead of directly appending the KKT conditions into upper-level constraints, we move part of the high-dimension constraints into the objective function and formulate them as penalty terms. The new penalty added objective and the adjusted constraints are written as Eq.~(\ref{eq:penalty_function}).

\begin{subequations}\label{eq:penalty_function}
\begin{align}
    \min_{\mathbf{p,x,v, \mu, \lambda, \beta, \pi}} \quad &\Phi'_{0} = -\sum_{s} \sum_{l} d_l p_l \, x_{s,l} \, q_{s} + \rho \sum_{l, s}\Lambda_{l,s} \label{eq:penalty_obj} \\
    \text{s.t.} \quad & \Lambda_{l,s} \geq 0, \quad \forall l\in A, s\in S \label{eq:relax_link_kkt}\\
    & \text{Eqs.~(\ref{eq:price_cap_single}-\ref{eq:sub_cap_single},\ref{eq:node_kkt}-\ref{eq:nonneg_single})}
\end{align}
\end{subequations}

The penalty term $\Lambda_{l,s}$ is the augmented first order condition Eq.~(\ref{eq:link_kkt}) for each ${l,s}$ pair. The penalty coefficient $\rho$ is a tunable hyperparameter that modifies the magnitude of the penalty. The higher the value of $\rho$, the heavier the penalty when not solved towards the lower-level optima. When $\rho$ is sufficiently high, Eq.~(\ref{eq:penalty_function}) is equivalent to minimizing the upper level objective while obtaining the lower level optimum. Therefore, we can also solve Eq.~(\ref{eq:penalty_function}) to obtain a global optimum of this bilevel problem.

\begin{proposition}
    There exists a value of $\rho$ for which the solution obtained from solving Eq.~(\ref{eq:penalty_function}) converges to the optimal solution obtained from Eq.~(\ref{eq:single_level}).
\end{proposition}

\begin{proof}
    The equality constraint Eq.~(\ref{eq:link_kkt}) can be written as Eq.~(\ref{eq:kkt_link_rewrite}).
\begin{subequations}\label{eq:kkt_link_rewrite}
\begin{align}
    & \Lambda_{l,s} \geq 0, \quad \forall l\in A, s\in S \label{eq:link_kkt_great}\\
    & \Lambda_{l,s} \leq 0, \quad \forall l\in A, s\in S \label{eq:link_kkt_less}
\end{align}
\end{subequations}
We denote the Lagrangian multiplier of Eq.~(\ref{eq:link_kkt_less}) to be $\eta_{l,s}$. We switch the original objective Eq.~(\ref{eq:single_obj}) to a minimization problem. We then relax Eq.~(\ref{eq:link_kkt_less}) to have Eq.~(\ref{eq:less_relax}).
\begin{subequations}\label{eq:less_relax}
\begin{align}
    \min \quad & - \sum_{s} \sum_{l} d_l p_l \, x_{s,l} \, q_{s} + \sum_{l,s} \eta_{l,s}\Lambda_{l,s} \\
    & \Lambda_{l,s} \geq 0, \quad \forall l\in A, s\in S \label{eq:link_kkt_great}
\end{align}
\end{subequations}
When $\rho \geq max(\eta)$, the penalty term $\rho \sum \Lambda$ dominates the incentive to violate Eq.~(\ref{eq:link_kkt_less}). Under this condition, any optimal solution of the original problem is optimal for the penalized problem. The equality constraint Eq.~(\ref{eq:link_kkt}) is recovered at the optimum even though the feasible set only enforces Eq.~(\ref{eq:link_kkt_great}). This completes the proof.
\end{proof}

To further reduce the computational load, we propose another algorithmic improvement. Since the lower-level problem is very efficient to solve due to its convexity, it is easy to obtain the lower-level optimum when \textbf{p} and \textbf{r} are fixed. We use an iterative update process to exploit this property and reduce the solution time for Eq.~(\ref{eq:penalty_function}). The overall proposed algorithm is summarized in \textbf{Algorithm 1}. The hyperparameters $\rho^0$, $\psi^-$, and $\psi^+$ control the rate of penalty update, where $\psi^+ > 1$ and $\psi^- < 1$. By starting with a small value $\rho^0$, Eq.~(\ref{eq:penalty_function}) is easier to solve to optimality. This comes with the cost of violating lower-level KKT condition by a large margin ($\sum \Lambda >>\epsilon$). By gradually increasing $\rho$ with $\psi^+$, the KKT violation shrinks, leading to the global optimum. However, when $\rho$ becomes too large, the single-level solution process can become computationally intensive. Therefore, we reduce $\rho$ by $\psi^-$ to control the step size. Within each iteration, the previously solved $\mathbf{p}$ and $\mathbf{r}$ are used to solve the lower-level problem first. This solution, along with all dual values and their bounds, can warm-start the single-level solution algorithm in a tighter space, leading to a shorter runtime. After meeting the convergence criteria or reaching the maximum number of iterations, we output the final values of $(\mathbf{p}^*, \mathbf{r}^*, \mathbf{x}^*, \mathbf{v}^*)$.

\begin{algorithm}
\caption{Penalty-Based Iterative Solution Method for Mobility Hub Platform Design Problem}
\begin{algorithmic}[1]
\State \textbf{Input:} Initial upper-level variable $\mathbf{p}^0, \mathbf{r}^0$, initial penalty parameter $\rho^0$, bound relaxation coefficient $\zeta$, penalty parameter update coefficient $\psi^-, \psi^+$, optimality gap target $\tau$, objective tolerance $\epsilon$, time limit $T$, maximum iterations $K$
\State $k \gets 0$
\Repeat
    \State \textbf{Step 1:} For given $\mathbf{p}^k, \mathbf{r}^k$, solve the lower-level problem Eq.~(\ref{eq:service_model}) to obtain optimal $\mathbf{x}^k, \mathbf{v}^k$
    \State \textbf{Step 2:} Obtain all Lagrangians $\mathbf{\lambda^k, \mu^k,\beta^k,\pi^k}$ from the lower level solutions and their current bounds.
    \State \textbf{Step 3:} Load and warm start the solution process for Eq.~(\ref{eq:penalty_function}) using $\mathbf{x^k, v^k, \lambda^k, \mu^k,\beta^k,\pi^k}$ and $\rho^k$. Update the bounds for $\mathbf{\lambda^k, \mu^k,\beta^k,\pi^k}$ by $\zeta$. Solve until time limit $T$ is reached or the optimality gap is lower than $\tau$. Obtain the new solution $\mathbf{p}, \mathbf{r}$.
    \State \textbf{Step 4:} Check optimality condition:
        \If{Optimality gap $ \leq \tau$}
            \State Go to step 5
        \Else
            \State $\rho^{k} \gets \psi^-\rho^{k}$
        \EndIf
    \State \textbf{Step 5:} Check penalty condition:
        \If{$\sum \Lambda \leq \epsilon$ \textbf{or} $k \geq K$}
            \State \textbf{Break}
        \Else
            \State $\rho^{k} \gets \psi^+\rho^{k}$
        \EndIf
    \State \textbf{Step 6:} $\textbf{p}^k, \textbf{r}^k=\textbf{p}, \textbf{r}$
    \State $k \gets k+1$
\Until{convergence}
\State Solve lower level problem Eq.~(\ref{eq:service_model}) with $\textbf{p}^*, \textbf{r}^*$. Obtain $\mathbf{x}^*, \mathbf{v}^*$
\State \textbf{Output:} Optimal solution $(\mathbf{p}^*, \mathbf{r}^*, \mathbf{x}^*, \mathbf{v}^*)$.
\end{algorithmic}
\end{algorithm}

\section{Numerical experiments}

In this section, we provide two sets of experiments to illustrate the formulation and the solution algorithm. We first use a toy network to verify the algorithm and illustrate the model. The second experiment is tested on a hypothetical multimodal network centered around three LIRR stations, which serve as potential mobility hubs for travelers commuting to New York City from Long Island.

\subsection{Illustrative example}
\begin{figure}
    \centering
    \includegraphics[width=0.75\linewidth]{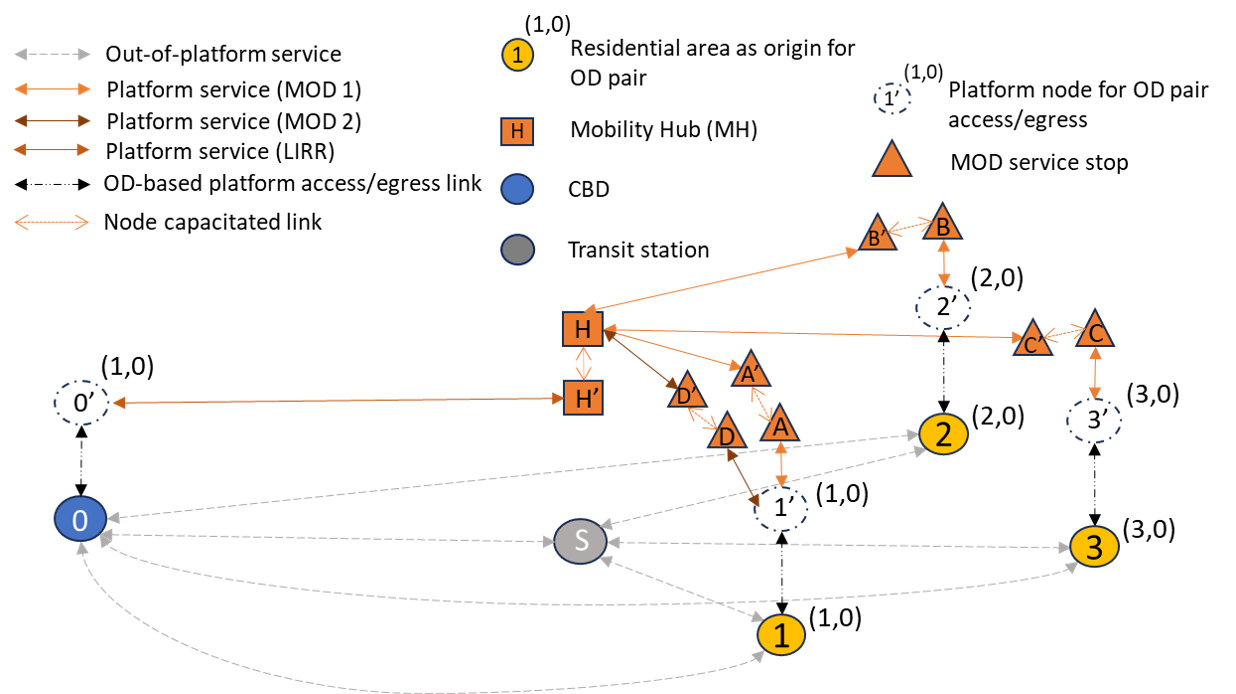}
    \caption{Toy network illustration}
    \label{fig:toy_network}
\end{figure}

A toy network is illustrated in \textbf{Fig.~\ref{fig:toy_network}}.Three OD pairs are considered: (1, 0), (2, 0), and (3, 0), each with a demand of 100 units over a typical operating period. The MaaS platform is considered for the transit line from mobility hub H to node $0'$, as well as MOD service links from nodes 1, 2, and 3 to node H. For node 1, two MOD service options are provided (e.g., car sharing and micromobility) with different cost attributes for both travelers and operators. The gray dotted links represent out-of-platform services (e.g., privately owned vehicles, walking, or other pay-as-you-go services). All platform-based nodes A, B, C, D, and H have virtual links to represent node capacities. Nodes A, B, and C are operated by the car sharing company (denoted as MOD 1) and node D is operated by the micromobility company (denoted as MOD 2). Both the MH node and MOD service nodes have capacity constraints and associated capacity costs. The MH capacity reflects the supply-side resources for supporting a fleet there: dedicated space for temporary stationing, maintenance, fueling, etc. The service node capacity reflects the availability constraint and the associated service stop establishment cost. For each OD pair, dedicated platform access links are defined for traveler pricing purposes.

Core input parameters are listed in \textbf{Table~\ref{tab: inputs}}. All access, egress, and virtual capacity links are assigned a small length value, which is used to prevent infeasible solutions without impacting the flow assignment results. Virtual node capacity links have zero cost and price and are also assigned a small length value to prevent infeasibility. The three MOD access links ($1-1'$, $2-2'$, $3-3'$) are controlled by the platform to adjust suitable service pricing in the upper-level problem along with operator-facing subsidies assuring operator profitability. Other services have fixed costs. $c_l^t$, $c_l^o$, and $p_l$ are all measured in monetary values. We want to highlight that after paying the access fee when traversing link $1-1'$, travelers from node 1 are free to choose either $1'-A$ or $1'-D$ without additional fare cost, which reflects the ''bundled'' feature of such a platform. Such behavioral decisions are determined by the lower-level problem, which decides the joint choice of traveler flow and hub capacities.

The toy example is solved on a device with an Intel(R) i7-13705 processor. We directly solve it to global optimality using Gurobi 12.1 and the branch-and-bound algorithm, with the optimality gap being less than $10^{-4}$. All scenarios reach the global optimum under 1 second.

% \textbf{BENCHMARK SCENARIO}: To better understand how the price allocation scheme impact both the MH based service flow and platform revenue, we run a benchmark scenario representing a traditional transit hub where operators may converge for transfers, but no digitalization is setup for cost transfers between the operators. This solution is obtained by simply solving the lower-level model only with a fixed price of $\$3$. When MOD service link prices are all fixed to \$3 per unit distance, only link (A, H) is actively used with a flow of 24.47 from node 1. Both nodes 2 and 3 are not served by the MOD service. Total revenue is \$73.4 and the proportion of opened MH capacity is 12.24\%.

The results are summarized in \textbf{Table~\ref{tab: toy_example_result}} under the columns "Price", "Subsidy", and "Flow". Flows on virtual capacity links (e.g., D-D') are not shown in the table for simplicity.
% \begin{remark}
% \textit{The model solution outputs  that is sensitive to link-level travel disutilities and operator costs while determining the.}
% \end{remark}
OD pair (3, 0) relies exclusively on the direct non-MaaS link, whereas OD pairs (1, 0) and (2, 0) exhibit multimodal routing behavior, combining both direct and MaaS-integrated links. Specifically, flows for OD (1, 0) and (2, 0) are distributed across three routes: the direct link, a non-MaaS path using the transit connection, and a MOD-only path. This highlights the benefits to travelers starting from nodes 1 and 2 when MOD services are introduced. However, under the provided pricing scheme for OD (1, 0), no travelers select service node A, and the service link chain $A-A'-H$ remains closed. This reflects the competitive dynamics of the market, where multiple services may be available but not all are viable for platform operation.

The platform access fees for OD (1, 0) and (2, 0) are set at 8.49 and 6.84, respectively, to maximize total revenue. To incentivize operators to participate in the platform without incurring losses, the operator receives a subsidy of 3.50 per traveler entering the service at node D for OD (1, 0). For OD pair (2, 0), the corresponding operator subsidy is 4.50. This indicates that involving operators in the route from service node B to hub H is more costly for the platform. For nodes A and C, other modes dominate the options. Therefore, no subsidy is offered, as opening services at these nodes is deemed unprofitable.

For OD (3, 0), the direct link in the dummy subnetwork dominates all other options. Even if the platform were to offer free access to the MOD service link (C,H), travelers would not choose it, nor would the operator serve that area. Due to the overall low usage of MaaS services in this scenario, the hub is significantly underutilized, operating at only 21.59\% of its maximum capacity. The total platform service profit in this case is \$161.10.

% \begin{remark}
% \textit{Compared to a benchmark transit hub, the subsidy-enabled MH can produce a higher capacity and increased revenues; the value of having a such a MH can be quantified as the difference in the objective value.}
% \end{remark}

% \textbf{ALTERNATIVE SCENARIO}: Now consider a scenario where the service price is capped to $\$2$ per unit distance. As shown in the "Price 2" and "Flow 2" columns of \textbf{Table~\ref{tab: toy_example_result}}, more travelers use the MOD service links on link (A, H) due to the lowered cap. For link (B, h), the price level remains the same, although with the price cap the required subsidy is reduced to only \$1.06. Again, no travelers choose link (C, H) due to the dominance of the direct link (3, 0). As a result, the platform sets the capacity to 20.31\% of maximum hub capacity to serve MaaS users. The MaaS platform earns a total of \$80. The \$10 dollar loss in revenue due to further constraining the maximum price translates to a higher enrollment in users coming from zone 1 at an overall lower social surplus. 

% \begin{remark}
% \textit{The model captures the social surplus and ridership impacts of a price cap.}
% \end{remark}

\begin{table}
\centering
\caption{Network Input and Parameters}
\label{tab: inputs}
\small % Reduces font size slightly to fit better
\setlength{\tabcolsep}{4pt} % Reduces white space between columns
\begin{tabular}{cccccc}
\toprule
\textbf{Source} & \textbf{Sink} & $\mathbf{c_l^t (\$)}$ & $\mathbf{c_l^o (\$)}$ & $\mathbf{\hat{p}_l (\$)}$ & $\mathbf{d_l}$ \\
\midrule
1   & 0   & 9 & 0 & N/A & 4   \\
2   & 0   & 9 & 0 & N/A & 4   \\
3   & 0   & 9 & 0 & N/A & 4   \\
1   & $1'$   & 0 & 0 & 10 & 0.1 \\
2   & $2'$   & 0 & 0 & 10 & 0.1 \\
3   & $3'$   & 0 & 0 & 10 & 0.1 \\
$1'$ & A & 0 & 0 & N/A & 0.1 \\
$1'$ & D & 0 & 0 & N/A & 0.1 \\
$2'$ & B & 0 & 0 & N/A & 0.1 \\
$3'$ & C & 0 & 0 & N/A & 0.1 \\
$A'$   & H   & 5 & 2 & N/A & 1   \\
$B'$   & H   & 5 & 2 & N/A & 2   \\
$C'$   & H   & 5 & 2 & N/A & 3   \\
$D'$  & H & 3 & 3 & N/A & 1 \\
1   & S   & 7 & 0 & N/A & 1   \\
2   & S   & 7 & 0 & N/A & 2   \\
3   & S   & 7 & 0 & N/A & 3   \\
$H'$   & $0'$   & 6 & 0 & N/A & 4   \\
S   & 0   & 7 & 0 & N/A & 4   \\
\midrule
% The descriptions now span 5 columns, leaving the 6th column for the value
\multicolumn{5}{l}{\textbf{Maximum hub capacity, $b_l$}} & \textbf{200} \\
\multicolumn{5}{l}{\textbf{Maximum service node capacity, $z_l$}} & \textbf{50} \\
\multicolumn{5}{l}{\textbf{Maximum subsidy allowed per traveler, $\hat{r}_l$}} & \textbf{5} \\
\multicolumn{5}{l}{$\mathbf{\alpha_1}, \mathbf{\alpha_2}$} & \textbf{1, 0.5} \\
\multicolumn{5}{l}{\textbf{Capacity cost, $c_i$}} & \textbf{1} \\
\bottomrule
\end{tabular}
\end{table}

\begin{table}
\centering
\caption{Comparison of link prices, subsidies and flows across base and alternative scenarios}
\label{tab: toy_example_result}
\begin{tabular}{llccc}
\toprule
\textbf{Link} & \textbf{OD Pair $s$} & \textbf{Price, $p_l$} & \textbf{Subsidy, $r_l$} & \textbf{Flow, $x_{l,s}q_s$}\\
\midrule
(1, 0)   & (1, 0) & N/A       & N/A   & 32.87  \\
(2, 0)   & (2, 0) & N/A       & N/A   & 82.59  \\
(3, 0)   & (3, 0) & N/A       & N/A   & 100    \\
(1, 1')  & (1, 0) & 8.49      & N/A   & 30.82 \\
(1', D)  & (1, 0) & N/A       & 3.5   & 30.82  \\
(D', H) & (1, 0) & N/A    & N/A   & 30.82  \\
(H', 0)  & (1, 0) & N/A       & N/A   & 30.82  \\
(2, 2')  & (2, 0) & 6.84     & N/A   & 12.35  \\
(2, B)  & (2, 0) & N/A      & 4.5   & 12.35  \\
(B', H) & (2, 0) & N/A    & N/A   & 12.35 \\
(H', 0)  & (2, 0) & N/A       & N/A   & 12.35  \\
(1, S)  & (1, 0) & N/A       & N/A   & 36.30  \\
(S, 0)  & (1, 0) & N/A       & N/A   & 36.30  \\
(2, S)  & (2, 0) & N/A       & N/A   & 5.06  \\
(S, 0)  & (2, 0) & N/A       & N/A   & 5.06  \\
\midrule
\multicolumn{2}{l}{\textbf{Hub Capacity, $v_i$}} & \multicolumn{2}{c}{21.59\%} \\
\multicolumn{2}{l}{\textbf{Objective value, $\Phi_0$}} & \multicolumn{2}{c}{$\$161.10$} \\
\bottomrule
\end{tabular}
\end{table}

\subsection{LIRR mobility hubs case study}

\begin{figure}[htbp]
    \centering
    \includegraphics[width=0.9\linewidth]{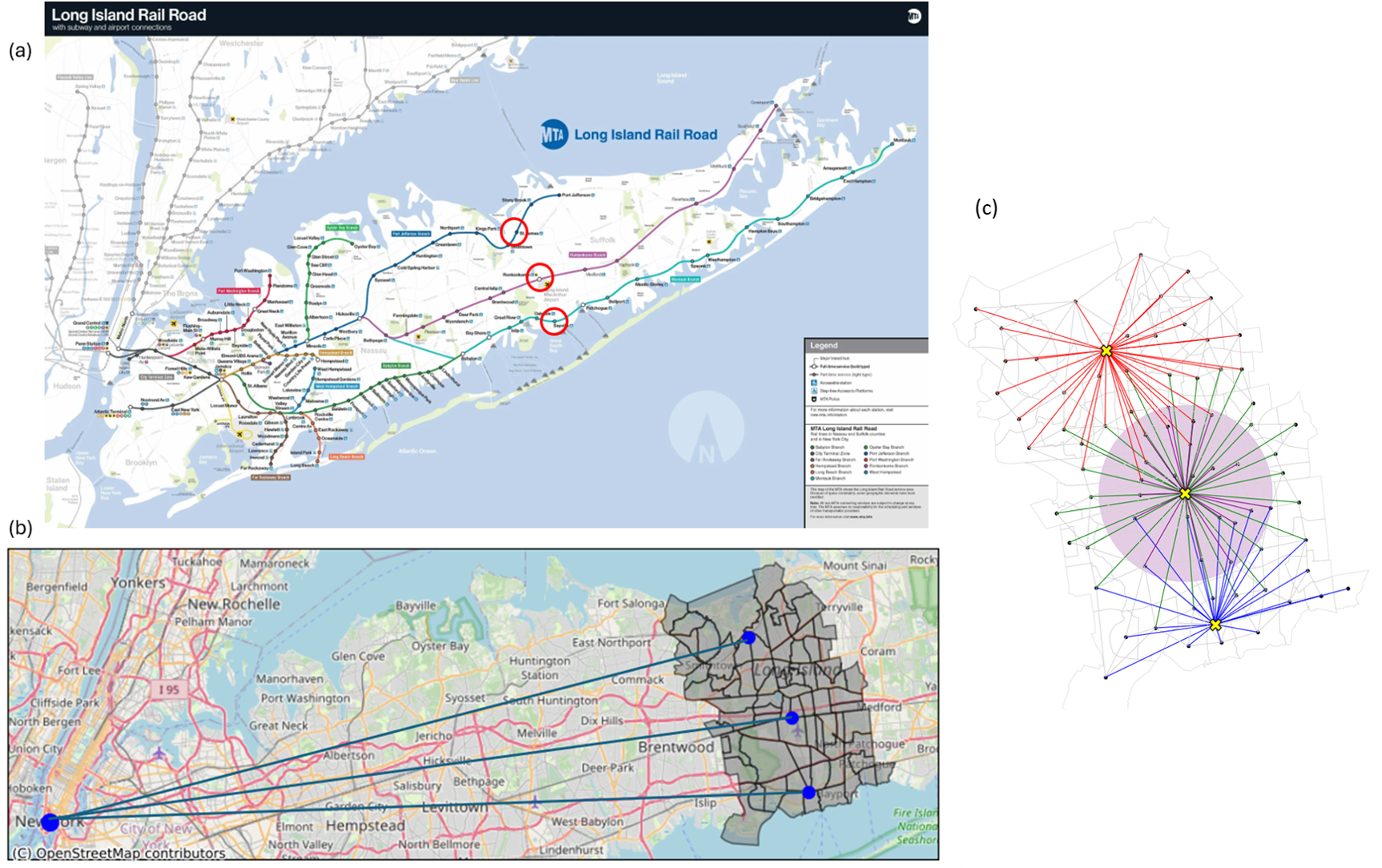}
    \caption{LIRR network illustration: (a) map of LIRR from the MTA with the three stations circled in red; (b) overlay of the census tracts within 5 miles of those stations; and (c) the MaaS subnetwork with links connecting each census tract to nearby stations}
    \label{fig:LIRR_network}
\end{figure}

We further use three LIRR stations: Ronkonkoma, St. James, and Sayville, and their surrounding neighborhoods, represented by travel analysis zones (TAZs), to test the proposed model. These three stations are treated as candidate mobility hubs. We use the centroids of census tracts within a 5-mile radius of each station as origin nodes, with Manhattan as the destination for a typical weekday morning AM period. In total, 78 neighborhoods are covered. Two layers of subnetworks are defined, similar to the toy network. \textbf{Fig.~\ref{fig:LIRR_network}} illustrates the platform subnetwork. Each service node is directly connected via MOD service links, with candidate mobility hubs established at each LIRR station. We assume these primary MOD services are managed by a single operator. To model market competition similar to the toy network, we introduce a second MOD service provider in the vicinity of Ronkonkoma station. This secondary operator functions within a 3-mile radius, covering 20 neighborhoods. In a real world context, the primary MOD service (MOD 1) represents larger scale microtransit options, such as ride hailing or ridepooling. For example, Suffolk County operates an on-demand service in Southampton that provides last mile access to the LIRR Southampton station and the S92 bus \citep{suffolktransit}. The localized MOD service (MOD 2) represents micromobility solutions, such as bike sharing or shared scooters. The resulting network comprises 98 service nodes connected by 133 direct microtransit service links. Additionally, three LIRR links are included to represent line-haul services. 

For the dummy subnetwork, 78 direct service links parallel to the MOD links are added to represent other access services. Another 78 direct links from service nodes to the Manhattan node are also included in the dummy subnetwork. As a result, the whole network has 380 nodes and 642 links serving the 78 OD pairs. By comparison, the expanded Sioux Falls network in \cite{liu2024demand} has 30 OD pairs, 82 nodes, and 748 links.

The 78 OD pairs, each originating from a service node to Manhattan, have demand simulated from a normal distribution with a mean of 60 and a standard deviation of 20, representing a typical weekday AM period (see \textbf{Table.~\ref{tab:simulated_demand}} in Appendix for details). The simulated total demand is 4,734, with an average being 60.69 per census tract OD pair. Note that this level of demand falls in the range of observation, as Ronkonkoma Station, the largest one, has 5,452 free parking spaces for general passengers, and fills up on a typical weekday. 

Each MOD service node is assigned a capacity of 50, while the Mobility Hub (MH) has a capacity of 2,000. The maximum platform access fee (which includes the LIRR ride) is capped at \$30, and the operator subsidy is limited to a maximum of \$5 per traveler. $\alpha_1$ and $\alpha_2$ are set to be 1 and 0.5, respectively. Other input elements are summarized in \textbf{Table~\ref{tab: LIRR_input}}. The interpretation of the $\hat{p_l}$ parameter varies depending on the link type. For MOD platform access links, it represents the price cap. Conversely, for dummy network links, it represents the fixed fee paid by travelers. $\bar{d}_l$ represents the average link length. To reflect the trade-off between micromobility and microtransit, MOD 2 is modeled with lower operating costs but slower service speeds than MOD 1.

For the experimental design, we first establish two baseline scenarios: one without the proposed platform ("business as usual") by solving the lower level problem without the platform subnetwork, and a second with the proposed platform by solving the problem with the stated network and input variables to determine the baseline objective value. We then perform a sensitivity analysis focusing on MOD operating costs within the platform network and driving costs within the dummy network to illustrate the sensitivities captured by the model. To further evaluate system dynamics, we then analyze three targeted scenarios: (1) eliminating micromobility links to assess operator participation effects; (2) removing the Ronkonkoma Mobility Hub designation to quantify its added value; and (3) reducing the operator subsidy cap to \$3 to examine the impact of tighter fiscal constraints on planning decisions.

All cases are run on a device with an Intel(R) i7-13705 processor. We initiate the proposed algorithm with $\rho_0$ equal to 300 and terminate the algorithm either after 10 iterations or when the optimality gap is below 0.1\%. Each iteration has a runtime limit of 10 minutes.

\begin{table}[htbp]
\centering
\caption{Subnetwork Mode Parameters}
\label{tab: LIRR_input}
\begin{tabular}{llccccc}
\toprule
\textbf{Subnetwork} & \textbf{Mode} & $\mathbf{c^t_l}(\$)$ & $\mathbf{c^o_l}(\$)$ & $\mathbf{\hat{p}_l}(\$)$ & $\mathbf{\bar{d}_l}(miles)$ \\
\midrule
\multirow{3}{*}{MOD}   & Access      & 0   & 0  & 20   & 0.1 \\
                       & MOD 1        & 0.42 & 0.5 & 0 & 3  \\
                       & MOD 2       & 0.5 & 0.3 & 0 & 1.8 \\
                       & LIRR        & 1.2 & 0   & 0 & 25 \\
\midrule
\multirow{3}{*}{Dummy} & LIRR dummy  & 1.30 & 0   & 20 & 25 \\
                       & Drive       & 1.33 & 0   & 21 & 30 \\
                       & Park-n-Ride & 0.5 & 0   & 0 & 5  \\
\bottomrule
\end{tabular}
\end{table}

\subsubsection{Base case}
Prior to solving the full bilevel model, we solved the lower-level model exclusively on the dummy network to establish traveler behavior in the absence of platform participation. Using the parameters listed in \textbf{Table~\ref{tab: LIRR_input}}, we observe a mode split where 64.52\% of travelers utilize the existing transit network and 35.45\% opt for direct driving. The lower level objective value is 5,451.97. These results serve as a "business as usual" baseline for comparison with the proposed bilevel model, allowing us to quantify the impact of introducing the platform service.

We then run the full bilevel model with the platform network added. The proposed algorithm finds a solution with an optimality gap of less than 1\% in $3092$ seconds. The runtime is significantly shorter than other link-based MaaS design models proposed by previous studies \citep{liu2024demand, yao2024design} when solving problems with similar scale. In comparison, the Sioux-Fall network-based case described in \cite{liu2024demand} requires longer than a 4-hour run time to converge. This illustrates the scalability of the proposed algorithm in real-world applications. 

\begin{figure}[htbp]
    \centering
    \includegraphics[width=0.75\linewidth]{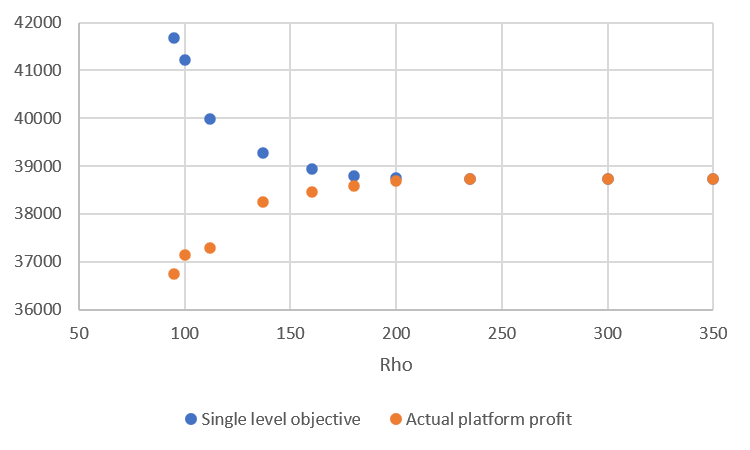}
    \caption{Convergence between single level objective and true platform profit.}
    \label{fig:convergence}
\end{figure}

To demonstrate convergence to the global optimum, we solved the model across a range of $\rho$ values. We then fix the upper-level decision variables to their corresponding optimal outputs and solved the lower-level model directly. Using the resulting flows and capacities, we calculated the actual platform profit. As illustrated in \textbf{Fig.~\ref{fig:convergence}}, both values converge as $\rho$ increases. However, the rate of convergence diminishes as $\rho$ approaches a critical value. Once $\rho$ exceeds this threshold, the optimum of the single-level formulation effectively represents the true global optimum of the bilevel model. Future studies could focus on developing algorithms that efficiently identify this critical threshold to further accelerate the solution process.

\begin{figure}[htbp]
    \centering
    \includegraphics[width=1\linewidth]{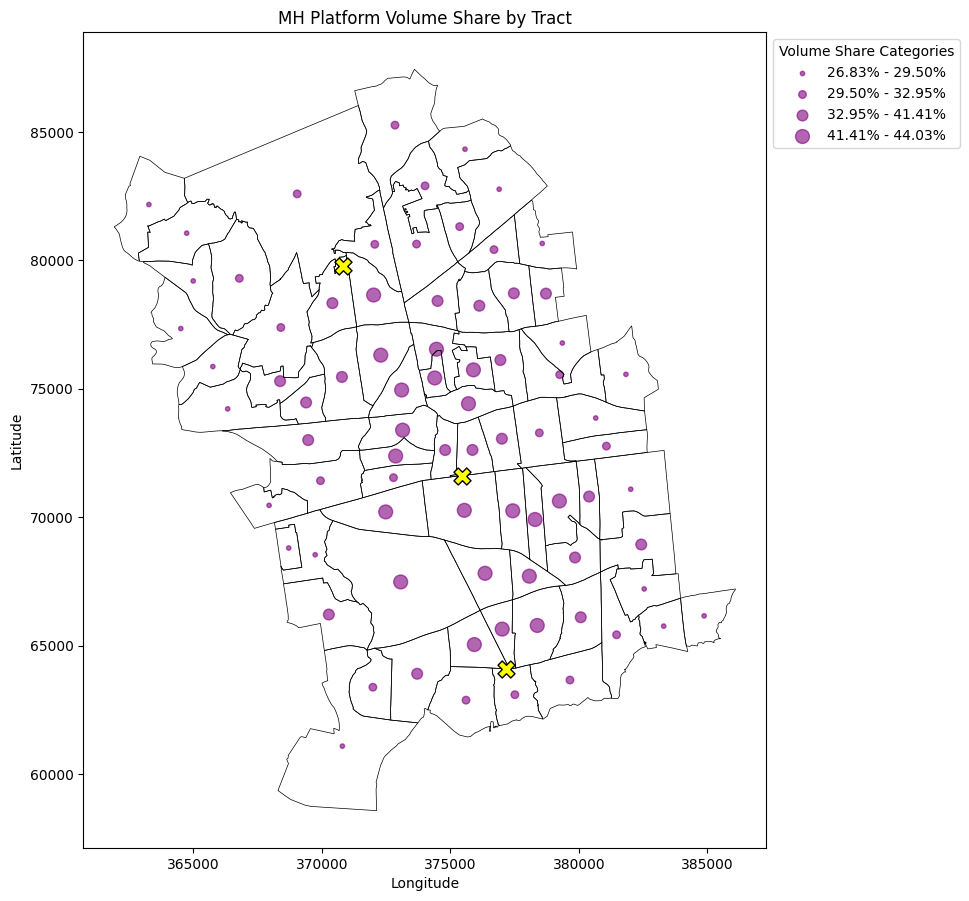}
    \caption{Mode share of MH based MOD services.}
    \label{fig:baselin_result_region}
\end{figure}

\begin{figure}[htbp]
    \centering
    \includegraphics[width=1\linewidth]{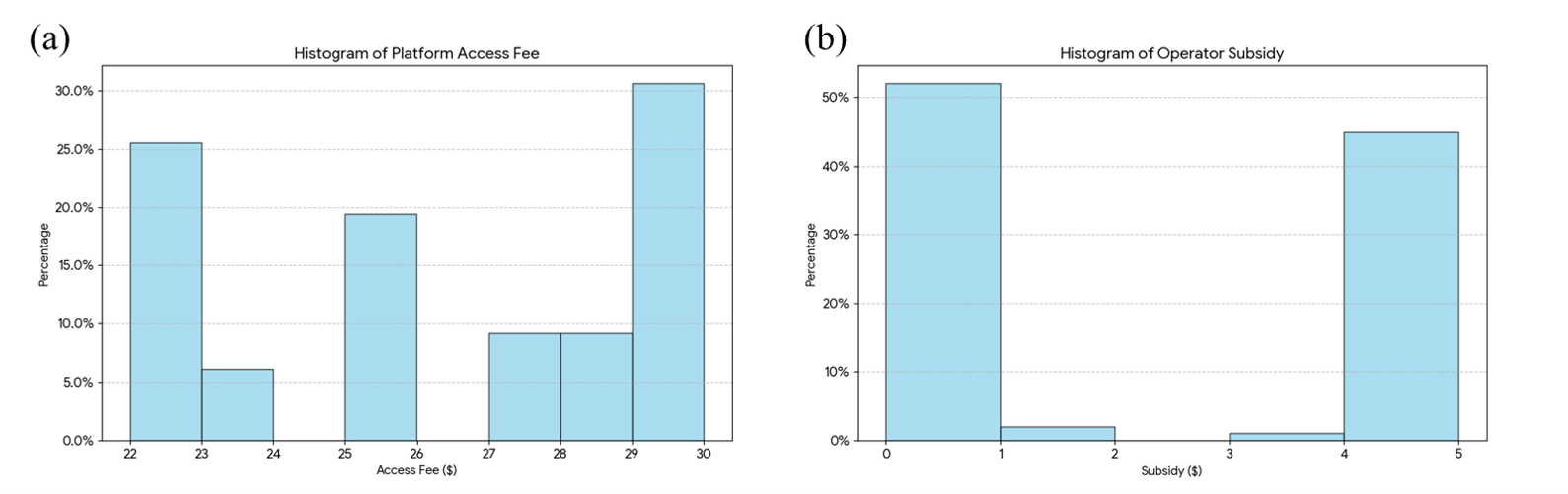}
    \caption{Histogram of (a) platform access fee and (b) histogram of operator subsidy.}
    \label{fig:baselin_result}
\end{figure}

The final solution has platform access prices (which includes the access mode and LIRR) ranging from \$22.47 to \$30, and the general distribution can be observed in \textbf{Fig.~\ref{fig:baselin_result}(a)}. This pricing structure is economically reasonable when compared to general LIRR commuting costs, ensuring the platform remains a viable option for travelers. The total profit from the MOD service is \$38721.23, and 35.32\% (1,661 passengers) of the total OD demand uses the MH centered platform. 741 people use service based at Ronkonkoma station, 550 use service based at St. James station, and 380 use service based at Sayville. \textbf{Fig.~\ref{fig:baselin_result_region}} illustrates the percent mode share of the MOD services across the 78 neighborhoods. 

The operator subsidy ranges from \$0 to \$5. However, the distribution is polarized, with most values falling either below \$1 or above \$4, as illustrated in \textbf{Fig.~\ref{fig:baselin_result}(b)}.This indicates that the platform strategically concentrates financial incentives on specific, high-priority service links rather than distributing them uniformly. We observe a strong positive correlation of 0.66 between operator subsidies and passenger flow. This suggests that the platform allocates higher per-traveler subsidies to high-demand links to ensure that the total subsidy volume is sufficient to cover the operator's aggregate costs. Conversely, the correlation between link length and subsidy level is weakly positive, at 0.11. While areas further from MHs generally require higher subsidies to be viable, distance alone is not the sole determinant. Areas easily accessible via the dummy network often have lower platform usage. Furthermore, variations in access fees create different passenger flows for areas at similar distances. Consequently, the relationship between distance and subsidy is only weakly correlated.

Both the major operator (MOD 1 in \textbf{Table~\ref{tab: inputs}}) and the smaller operator receive the subsidies to reach break even. The major operator incurs \$4781.41 in operating costs for the in-platform service and receives an equivalent \$4781.41 in subsidies, and the smaller operator incurs \$220.76 in operating costs and receives \$220.76 from the platform. With the introduction of the platform service, the mode share of direct driving drops to 22.46\%, representing a 36.6\% decrease compared to the no-platform scenario. To quantify the utility impact, we compared the lower-level objective values of both models. The objective value decreases from 5,451.97 in the baseline (dummy network only) scenario to 5,154.90 in the optimal platform case. Since the lower-level objective measures user disutility, this reduction demonstrates the platform's positive impact on overall system welfare. Per Remark 1, this particular baseline suggests there is an opportunity for a subsidy-enabled MH platform to serve this market. 

\begin{remark} \label{remark: zero_subsidy}
\textit{The model identifies an opportunity for a MH platform when it reduces the overall travel disutilities with a fare bundle that can cover the costs of all participating operators.}
\end{remark}

We subsequently conducted two sensitivity analyses focusing on cost-related components: the MOD link operating cost and the direct driving cost. These analyses reveal how fluctuations in these costs impact both the platform's profitability and the utilization of the three MHs.

\paragraph{MOD link operating cost sensitivity analysis} \mbox{} \\

\begin{figure}
    \centering
    \includegraphics[width=0.75\linewidth]{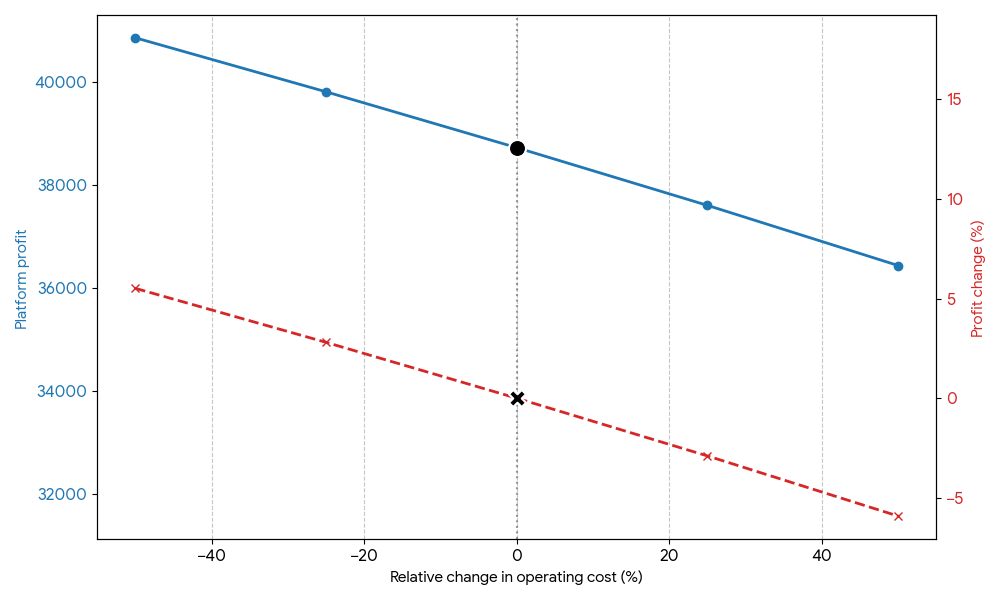}
    \caption{Relationship between platform profit and MOD link operating cost}
    \label{fig:operating_cost_sensitivity}
\end{figure}

\noindent We first perform a sensitivity analysis on the operating costs of MOD links. We varied the operating costs for services operated by both MOD 1 and MOD 2 across a range of -50\% to +50\%. \textbf{Fig.~\ref{fig:operating_cost_sensitivity}} illustrates the relationship between these operating costs and the platform profit. As operating costs rise, the platform is required to provide higher subsidies. This leads to an increase in the platform access fee, as shown in \textbf{Table.~\ref{tab:cost_impact}}, which subsequently reduces platform service usage (and vice versa). Notably, we observed a linear relationship between operating costs and platform profit. A 50\% fluctuation in operating cost results in a 5\% shift in platform profit under these baseline settings. Given that the subsidy comprises only 13\% of the total platform revenue in the base scenario, even a 50\% variation in cost leads to only a marginal adjustment in the platform access fee.

\begin{table}
\caption{Impact of operating cost changes on access fee and share of platform usage}
\label{tab:cost_impact}
\centering
% \resizebox{width}{height}{content}
% We use \columnwidth for width and ! for height to keep aspect ratio
\resizebox{\columnwidth}{!}{%
    \begin{tabular}{cccc}
    \toprule
    \textbf{Operating Cost Change (\%)} & \textbf{Average Access Fee} & \textbf{Relative Fee Change (\%)} & \textbf{Share (\%) of MH usage} \\
    \midrule
    -50 & 26.13 & -1.94 & 36.14 \\
    -25 & 26.37 & -1.06 & 35.79 \\
    0     & 26.65 & 0.00  & 35.32 \\
    +25 & 26.92 & 1.01  & 34.82 \\
    +50 & 27.16 & 1.93  & 34.30 \\
    \bottomrule
    \end{tabular}%
}
\end{table}

\begin{figure}
    \centering
    \includegraphics[width=0.75\linewidth]{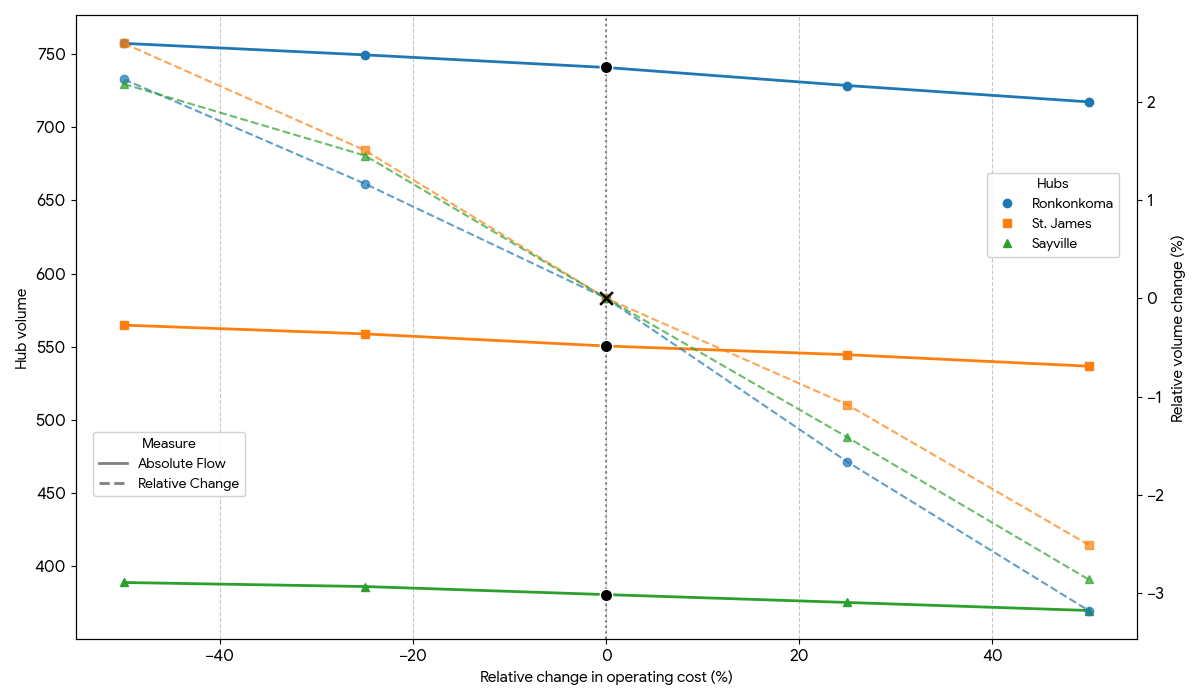}
    \caption{Relationship between MH usage and MOD link operating cost}
    \label{fig:hub_volumn_operating_cost}
\end{figure}

\textbf{Fig.~\ref{fig:hub_volumn_operating_cost}} illustrates the relationship between link operating costs and the volumes across the three MH locations. While the volumes at all three hubs exhibit a consistent inverse relationship with operating costs, the response is asymmetric. The negative impact of a cost increase is stronger than the positive gain from a cost reduction. Specifically, a 50\% decrease in operating cost increases hub volumes by 2\% to 2.5\%, with St. James benefiting the most. However, a 50\% cost surge leads to a sharper decline of 3.0\%, with Ronkonkoma becoming the most impacted. This suggests that while lower costs encourage some additional usage, higher costs penalize demand more severely.

\paragraph{Direct driving cost sensitivity analysis} \mbox{} \\

\begin{figure}
    \centering
    \includegraphics[width=0.75\linewidth]{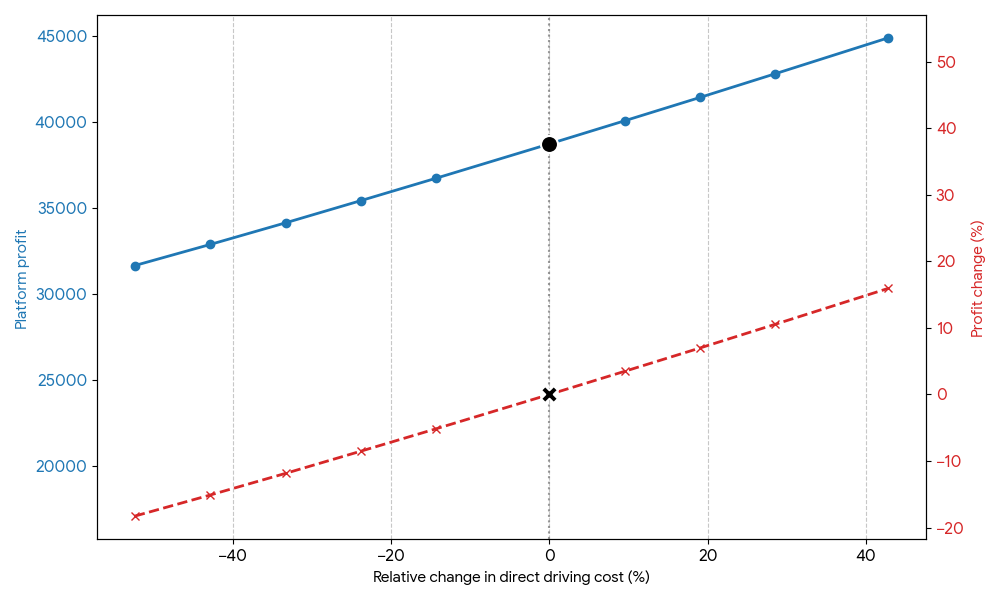}
    \caption{Relationship between platform profit and direct driving cost}
    \label{fig:driving_cost_sensitivity}
\end{figure}

\noindent We subsequently perform a sensitivity analysis on direct driving costs to evaluate the impact of competing mode pricing. Fig.~\ref{fig:driving_cost_sensitivity} illustrates the resulting relationship. Since private driving functions as a substitute mode, increasing its cost dissuades travelers from using private vehicles. This drives higher demand for MH services and increases platform profit. For instance, a 43\% increase in driving costs, which is comparable to a \$9 congestion toll for entering New York City, results in a 16\% gain in platform profit. Consistent with our operating cost analysis, the relationship between platform profit and driving cost changes is linear. Table.~\ref{tab:driving_cost_impact} details the corresponding shifts in access fees and the platform's overall mode share. Notably, the platform access fee covaries positively with driving costs. As driving becomes more expensive, the platform strategically raises its access fee to capitalize on the influx of deterred drivers. This creates an indirect effect where travel becomes universally more expensive, negatively impacting travelers who rely on the platform as an alternative to driving.

\begin{table}[htbp]
\caption{Impact of driving cost changes on access fee and platform mode share}
\label{tab:driving_cost_impact}
\centering
\small % Slightly smaller font often helps in academic papers
\begin{tabular}{cccc}
\toprule
\makecell{\textbf{Driving Cost} \\ \textbf{Change (\%)}} & 
\makecell{\textbf{Avg. Access} \\ \textbf{Fee}} & 
\makecell{\textbf{Relative Fee} \\ \textbf{Change (\%)}} & 
\makecell{\textbf{Mode Share} \\ \textbf{(\%)}} \\
\midrule
-52 & 24.59 & -7.72 & 31.66 \\
-43 & 25.00 & -6.17 & 32.27 \\
-33 & 25.42 & -4.63 & 32.88 \\
-24 & 25.83 & -3.08 & 33.48 \\
-14 & 26.19 & -1.71 & 34.17 \\
0   & 26.65 & 0.00  & 35.32 \\
10  & 26.93 & 1.06  & 36.10 \\
19  & 27.21 & 2.11  & 36.89 \\
29  & 27.49 & 3.16  & 37.68 \\
43  & 27.88 & 4.60  & 38.90 \\
\bottomrule
\end{tabular}
\end{table}

\begin{figure}[htbp]
    \centering
    \includegraphics[width=0.75\linewidth]{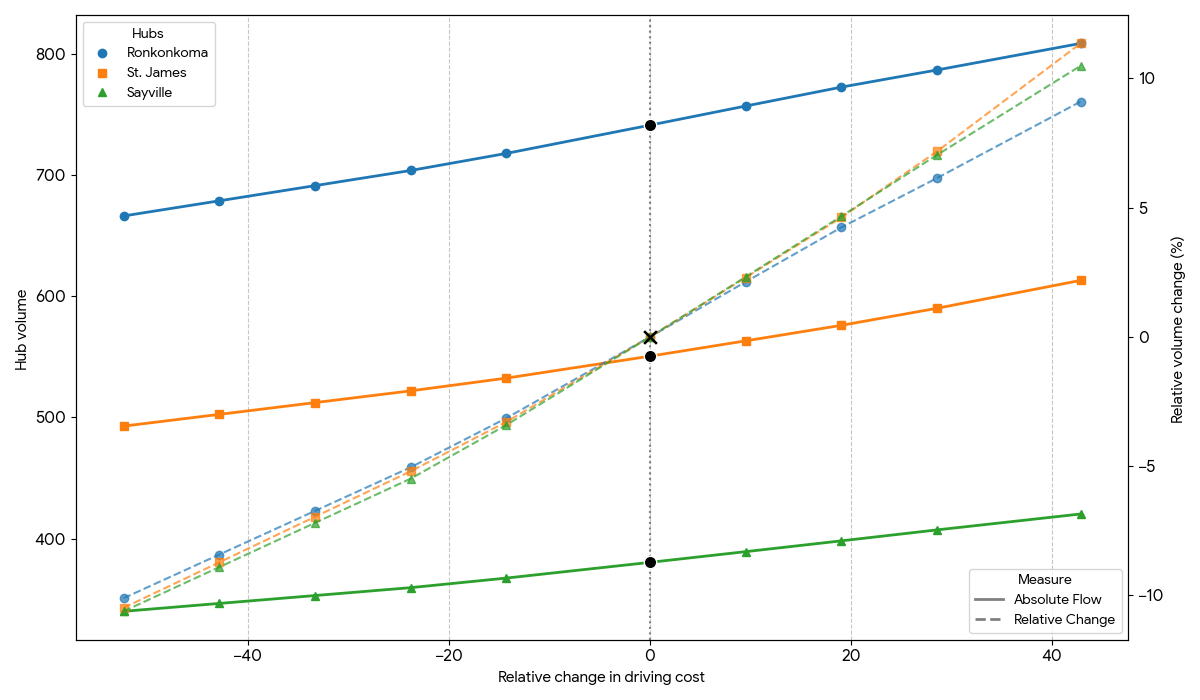}
    \caption{Relationship between MH usage and direct driving cost}
    \label{fig:hub_volume_driving_cost}
\end{figure}

Fig.~\ref{fig:hub_volume_driving_cost} illustrates the relationship between direct driving costs and traffic volumes across the three MH locations. When driving costs decrease, the relative decline in volume is nearly identical across all three hubs. However, as driving costs rise, the relative growth rates begin to diverge slightly. St. James exhibits the greatest relative gain, while Ronkonkoma shows the lowest. Despite this minor divergence, the three hubs demonstrate broadly similar demand elasticity with respect to driving costs. Consequently, policy measures that change the cost of driving would likely produce a consistent impact on ridership across the three hubs.

\subsubsection{Evaluating the impact of in-platform service competition}
To examine the impact of reduced operator competition, we eliminated MOD 2 service entirely. This consolidation resulted in a 0.73\% decline in total platform profit (\$38,439.05) and a marginal 0.1\% drop in ridership (1,659 passengers). While the average platform access fee remained relatively stable (increasing from \$26.65 to \$26.72), the pricing distribution shifted significantly. The first quartile rose from \$22.95 to \$25.29, and the median increased from \$25.86 to \$27.88. This indicates that reduced competition disproportionately raised costs for lower-priced trips. 

Due to the higher operating costs associated with the remaining operator, the total subsidy rises from \$5002.17 to \$5072.32, an increase of 1.40\%. Therefore, incorporating more service options into the platform proves beneficial for both the platform and travelers. The platform can potentially reduce total subsidy costs, while travelers benefit from greater access and a wider range of choices within the platform ecosystem with lower access fee.

\begin{remark}
\textit{Reduced competition in participating access operators to the MH platform can disproportionately raise costs for travelers even as the impact to the platform's profits may be minimal.}
\end{remark}

\subsubsection{Quantifying the value of an MH connecting microtransit to Ronkonkoma station}
To quantify the value added of the MH at Ronkonkoma, we apply the same model to a network modification where the Ronkonkoma MH (the central hub in \textbf{Fig.~\ref{fig:LIRR_network}}) and its associated MOD links are removed. Note that the Ronkonkoma station itself remains accessible within the dummy network. In this scenario, 63 MOD service nodes and 63 MOD links remain available on the MaaS platform, with only 15 census tracts losing platform service access.

Consequently, MaaS platform usage drops significantly from 35.32\% in the baseline to 22.30\%, with travelers diverting to the two remaining MHs or alternative out-of-platform modes. Total platform revenue contracts to \$24,311.86, implying that the inclusion of the Ronkonkoma MH generates an incremental value of \$14,409.37 per AM peak period. Despite these volume shifts, the pricing structure remains stable. Average access fees change by less than 1\%, with the first quartile shifting slightly from \$22.95 to \$23.38 and the median remaining virtually unchanged (\$25.86 vs. \$25.84).

\begin{remark} \label{remark: link_price_max}
\textit{Because the lower-level model is a coalitional choice model, the change in objective value for removing a MH captures the social surplus value of the MH that travelers and operators are collectively willing to pay for, i.e. the compensating variation of the MH location alternative.}
\end{remark}

\subsubsection{Lower subsidy cap provided by platform}
We reduce the subsidy cap from \$5 to \$3 while keeping the original network structure to evaluate the impact of subsidy constraints on platform planning decisions. The platform profit declines to \$37,746.32, representing a 2.52\% decrease compared to the baseline scenario, while the total number of travelers using the platform drops by 0.57\% to 1652. Although the average access fee remains relatively stable with marginally decreasing by 1.16\% to \$26.34, the allocation of subsidies shifts drastically. The total subsidy allocated to the two operators amounts to \$4,902.63 with high consolidation towards the major operator. MOD 1 sees a 1\% increase in subsidy allocation with \$4,829.37, while the smaller MOD 2 faces a significant reduction of 66.81\% to receiving only \$73.26. A sharp cut in subsidies caused MOD 2 ridership to fall from 149 to 42 and to only serve 15 out of the 20 covered areas, pushing the operator to the brink of market exit.

This drastic shift indicates that when subsidy flexibility is capped, the platform creates a "flight to quality", prioritizing the major operator (MOD 1) to secure core network capacity. The smaller operator (MOD 2), which provides a cheaper but slightly inferior service, likely relies on higher marginal incentives to keep operating. When the subsidy cap prevents the platform from bridging that gap, the platform effectively abandons the niche operator to preserve the financial viability of the primary fleet. Consequently, the restrictive cap reduces overall service diversity and accessibility, leading to a net loss in platform access and profitability. These results further reinforce the importance of model analysis of operating policies surrounding MHs, as the impacts may not be easily predicted.

\begin{remark}
\textit{Restrictive subsidy caps can impact different operators disproportionately.}
\end{remark}

\section{Conclusion}
In this study, we propose a MH platform design model addressing both the theoretical and computational gaps identified in recent literature for both MaaS and mobility hub based applications. While previous research has largely emphasized the role of MHs as physical transfer points or focused solely on optimal facility location, our approach extends the function of MHs to include their strategic influence on both operator and traveler decision-making, particularly through the explicit modeling of pricing structures and targeted subsidies. By leveraging the PURC framework within a bilevel optimization structure, we capture the complex interplay between public platforms and private mobility operators while also ensuring computational scalability necessary for real-world, large-scale applications.

The proposed model integrates the joint operational planning between travelers and operators, which is essential for effective MaaS platform deployment, while also capturing the competitive dynamics of multiple operator participation. The framework utilizes a bilevel structure: a revenue-maximizing upper level, where the platform determines optimal traveler-facing pricing and operator-facing subsidy allocation strategies. A convex quadratic lower level, where travelers and operators jointly optimize flows and capacities. We transform this bilevel formulation into a single-level problem using KKT conditions, guaranteeing the existence of a global optimum. Furthermore, an augmented penalty-based algorithm is proposed to enhance computational tractability and reduce runtime.

We demonstrate the model application through two sets of numerical examples. The toy example demonstrates the effectiveness of the proposed model and validates the existence of a global optimum. The results show that the model can optimally allocate link-specific resources from a revenue-maximization perspective. We further test the model on a larger network based on LIRR stations and their service areas. A total of 380 nodes and 642 links are involved in the base network with 78 OD pairs. Several insights are obtained through the tests:
\begin{itemize}
    \item The proposed model is much more scalable compared to similar platform design models in past studies. Similar-sized problems require hours of runtime (\citep{liu2024demand, yao2024design}), while ours takes less than 1 hour when setting the right conditions.
    \item Platform profit exhibits a linear relationship with both MOD link operating costs and direct driving costs within the analyzed network configuration.
    \item The platform's profitability demonstrates greater sensitivity to changes in direct driving costs than to changes in MOD service operating costs.
    \item All three MHs demonstrate similar level of ridership shifts when subjected to variations in service operating and driving costs.
    \item Market competition is important in reducing both the overall subsidy level and platform access pricing, specifically when the operators control different levels of market power.
    \item The use of the PURC-based lower-level problem allows policymakers to quantify the social surplus value of a MH, and to determine how to control subsidies to encourage platform service participation and traveler access.
\end{itemize}
The proposed model can be applied to more complex MaaS service designs involving more operators. With more operators, subsidy schemes based on setting a rate per operator can be done. Heterogeneous travel groups can be considered (e.g. splitting each OD pair into portions by income group). The proposed model is sensitive to changes in network topology, a common feature in existing service design problems. When different spatial aggregations are used to construct the network, associated parameters, including link travel times, operating costs, and coefficients such as $\alpha_1$ and $\alpha_2$, must be adjusted and recalibrated accordingly. Therefore, proper network construction using accurate observational data is an essential prerequisite. In future studies, we aim to incorporate supply-side functions that reflect operational performance into the bi-level model. This would relax the current assumption of fixed supply-side costs, allowing the model to determine optimal operator strategies corresponding to the platform's planning decisions.

There are multiple gaps remaining. Congestion effects on out-of-platform links can also be added. The link-additive approach could negate non link-additive factors impacting traveler and operator choices. Competition between platforms can be considered by changing the upper level into a generalized Nash equilibrium. Alternative applications of mobility hubs can be studied using this framework: urban air mobility and freight distribution, for example. The model can be further expanded to three-sided markets to deal with electric charging integration.

\section*{Acknowledgments}
The project is funded by the National Science Foundation (NSF) CMMI-2423908.

\bibliographystyle{apalike} 
\bibliography{reference}

\section*{Appendix}
\setcounter{table}{0}
\renewcommand{\thetable}{A\arabic{table}}
\begin{longtable}{l c c}
\caption{Simulated Demand Data} \label{tab:simulated_demand} \\
\toprule
Origin (Census Tract) & Destination & Demand \\
\midrule
\endfirsthead

\toprule
Origin (Census Tract) & Destination & Demand \\
\midrule
\endhead

\midrule
\multicolumn{3}{r}{\textit{Continued on next page}} \\
\endfoot

\bottomrule
\endlastfoot

3610300134704 & NYC & 37 \\
3610300134902 & NYC & 75 \\
3610300134903 & NYC & 50 \\
3610300134904 & NYC & 72 \\
3610300134906 & NYC & 44 \\
3610300134907 & NYC & 56 \\
3610300135002 & NYC & 76 \\
3610300135003 & NYC & 57 \\
3610300135004 & NYC & 65 \\
3610300135005 & NYC & 31 \\
3610300135104 & NYC & 60 \\
3610300135301 & NYC & 52 \\
3610300135303 & NYC & 87 \\
3610300135304 & NYC & 97 \\
3610300135401 & NYC & 46 \\
3610300135402 & NYC & 35 \\
3610300135403 & NYC & 83 \\
3610300145704 & NYC & 47 \\
3610300145803 & NYC & 41 \\
3610300145804 & NYC & 90 \\
3610300145805 & NYC & 64 \\
3610300145807 & NYC & 78 \\
3610300145808 & NYC & 52 \\
3610300146402 & NYC & 83 \\
3610300146403 & NYC & 83 \\
3610300146404 & NYC & 72 \\
3610300146604 & NYC & 53 \\
3610300146605 & NYC & 49 \\
3610300146606 & NYC & 58 \\
3610300146607 & NYC & 62 \\
3610300146608 & NYC & 40 \\
3610300146611 & NYC & 67 \\
3610300146612 & NYC & 42 \\
3610300146613 & NYC & 109 \\
3610300146614 & NYC & 67 \\
3610300146615 & NYC & 64 \\
3610300147503 & NYC & 72 \\
3610300147601 & NYC & 87 \\
3610300147602 & NYC & 58 \\
3610300147701 & NYC & 47 \\
3610300147702 & NYC & 28 \\
3610300147802 & NYC & 63 \\
3610300147803 & NYC & 47 \\
3610300147804 & NYC & 57 \\
3610300147901 & NYC & 58 \\
3610300147902 & NYC & 44 \\
3610300158002 & NYC & 76 \\
3610300158006 & NYC & 72 \\
3610300158007 & NYC & 29 \\
3610300158009 & NYC & 33 \\
3610300158010 & NYC & 51 \\
3610300158011 & NYC & 29 \\
3610300158103 & NYC & 72 \\
3610300158104 & NYC & 48 \\
3610300158107 & NYC & 77 \\
3610300158108 & NYC & 87 \\
3610300158110 & NYC & 97 \\
3610300158115 & NYC & 55 \\
3610300158502 & NYC & 58 \\
3610300158505 & NYC & 53 \\
3610300158506 & NYC & 78 \\
3610300158507 & NYC & 67 \\
3610300158508 & NYC & 62 \\
3610300158510 & NYC & 31 \\
3610300158511 & NYC & 65 \\
3610300158512 & NYC & 56 \\
3610300158604 & NYC & 78 \\
3610300158605 & NYC & 51 \\
3610300158606 & NYC & 62 \\
3610300158607 & NYC & 67 \\
3610300158608 & NYC & 74 \\
3610300158609 & NYC & 17 \\
3610300158802 & NYC & 75 \\
3610300158803 & NYC & 46 \\
3610300158804 & NYC & 64 \\
3610300158900 & NYC & 68 \\
3610300159000 & NYC & 59 \\
3610300159201 & NYC & 72 \\
\end{longtable}

\end{document}